\documentclass[reqno,11pt]{amsart}

\usepackage[top=2.0cm,bottom=2.0cm,left=3cm,right=3cm]{geometry}
\usepackage{amsthm,amsmath,amssymb}
\usepackage{mathrsfs,amsfonts,dsfont,functan,extarrows,mathtools}
\usepackage[colorlinks]{hyperref}

\usepackage{marginnote}
\usepackage{xcolor}

\newtheorem{theorem}{Theorem}[section]
\newtheorem{proposition}[theorem]{Proposition}

\newtheorem{lemma}[theorem]{Lemma}

\numberwithin{equation}{section}
\allowdisplaybreaks

\def\ls{\lesssim}

\def\ep{\varepsilon}

\def\p{\partial}
\def\div{{\rm div}}
\def\d{\mathrm{d}}
\def\no{\nonumber}
\def\L{\mathcal{L}}
\def\E{\mathcal{E}}
\def\SE{\sup_{s\in[0,t]}\mathcal{E}}
\def\wr{\widetilde{\rho}}
\def\wu{\widetilde{u}}
\def\wg{\widetilde{g}}
\def\WE{\widetilde{E}}
\def\WD{\widetilde{D}}
\def\WWE{\widetilde{\mathcal{E}}}

\newcounter{wronumber}

\newcommand{\abs}[1]{\left|#1\right|}
\newcommand{\nm}[1]{\left\|#1\right\|}

\newcommand{\lx}{{L^2_x}}
\newcommand{\lxq}{{L^2_{x,M\!q}}}
\newcommand{\hx}[1]{{H^ #1 _x}}
\newcommand{\hxq}[1]{{H^ #1 _x L^2_{M\!q}}}
\newcommand{\hxqm}[1]{{\mathcal{H}^ #1 _{x,M\!q}}}
\newcommand{\hxqe}[1]{{\mathcal{H}^{#1,\eta}_{x,M\!q}}}

\newcommand{\skp}[2]{\left( #1,\, #2 \right)}
\newcommand{\skt}[2]{( #1,\, #2 )}
\newcommand{\skpm}[2]{\left\langle #1,\, #2 \right\rangle_M}

\newcommand{\agl}[1]{\left\langle #1 \right\rangle}

\begin{document}

\title[Global Small Solutions of Compressible Polymers]{Global classical solutions to a Compressible Model for Micro-Macro Polymeric Fluids Near Equilibrium}

\author[N. Jiang]{Ning Jiang}
\address[Ning Jiang]{\newline School of Mathematics and Statistics, Wuhan University, Wuhan, 430072, P. R. China}
\email{njiang@whu.edu.cn}

\author[Y. Liu]{Yanan Liu}
\address[Yanan Liu]{\newline School of Statistics and Mathematics, Zhongnan University of Economics and Law, Wuhan, 430073, P. R. China}
\email{liuyn@znufe.edu.cn}

\author[T.-F. Zhang]{Teng-Fei Zhang}
\address[Teng-Fei Zhang]{\newline School of Mathematics and Physics, China University of Geosciences, Wuhan, 430074, P. R. China}
\email{zhtengf@mail2.sysu.edu.cn}

\date{}

%%%%%%%%%%%%%%%%%%%%%%%%%%%%%%%%%%%%%%%%%%%%%%%%%%%%%%%%%%%%%%%%%%%%%%%%%%%%%%
%%%%%%%%%%%%%%%%%%%%%%%%%%%%%%%%%%%%%%%%%%%%%%%%%%%%%%%%%%%%%%%%%%%%%%%%%%%%%%

\begin{abstract}
In this paper we first employ the energetic variational method to derive a micro-macro model for compressible polymeric fluids. This model is a coupling of isentropic compressible Navier-Stokes equations with a nonlinear Fokker-Planck equation. We then prove the global in time existence of the smooth solution near the global equilibrium.

\vspace*{5pt}
\noindent{\it Keywords}: Compressible polymeric model; Micro-macro system; Global existence; Classical solutions; Small initial data. 

\noindent{\it MSC 2010}: 35A01, 35Q30, 35Q84, 76N10

\end{abstract}

\maketitle

% \tableofcontents

%%%%%%%%%%%%%%%%%%%%%%%%%%%%%%%%%%%%%%%%%%%%%%%%%%%%%%%%%%%%%%%%%%%%%%%%%%%%%%
%%%%%%%%%%%%%%%%%%%%%%%%%%%%%%%%%%%%%%%%%%%%%%%%%%%%%%%%%%%%%%%%%%%%%%%%%%%%%%
\section{Introduction} % (fold)
\label{sec:intro}

% section introduction (end)

The research on the non-Newtonian (complex) fluids has been an active field in recent years, in which the viscoelastic fluids (such as polymeric fluids, liquid crystals and so on) are one of the most important types. These material possess elastic and viscous behavior, and fluid properties as well, hence they can be considered as the intermediate state between solid and fluid. The viscoelastic fluids are commonly modeled by systems coupling of fluids equations and kinetic equation, and hence involves the multi-scale properties.

In present paper we study a micro-macro model for compressible polymeric fluids, in which a polymer is viewed as an bead-spring model consisting of two beads connected by an elastic spring represented by a vector $q$ (the elongation of two beads), see \cite{BAH-1987-BOOK,BCAH-1987-BOOK}. The polymer particles are described by a probability density function $\Psi(t,x,q)$ satisfying that $\int_{D} \Psi(t,x,q) \d q =1$, which represents the distribution of particles' elongation vector $q \in D$ ($D$ is an bounded or unbounded domain in $\mathbb{R}^3$), at time $t\ge 0$ and position $x\in \Omega \subset \mathbb{R}^3$, with its evolution governed by a nonlinear Fokker-Planck equation. Specifically, we consider the case that the polymer particles are transported by a compressible fluid flow, which are governed by the following micro-macro coupling system:
	\begin{align}\label{eq:C-polymeric}
	  \begin{cases}
	  	\partial_t \varrho + \div (\varrho u) = 0,
		\\[5pt]
				\partial_t (\varrho u) + \div (\varrho u \otimes u) + \tfrac{1}{\text{Ma}^2} \nabla_x P(\varrho)
				= \div \Sigma(u) + \tfrac{1}{\text{De}} \tfrac{\lambda}{r} \div \int_{\mathbb{R}^3} (\nabla_q U \otimes q \Psi) \d q ,
	  \\[7pt]
	  	\partial_t \Psi + \div (u \Psi) + \nabla_q \cdot (\nabla_x u q\Psi)
		 	=\tfrac{1}{\text{De}} \nabla_q \cdot \left( \sigma \nabla_q \Psi + \tfrac{1}{r} \nabla_q U \Psi \right),
	  \end{cases}
	\end{align}
where $(\varrho,u,\Psi)= (\varrho(t,x),u(t,x),\Psi(t,x,q))$ are the density, the velocity field and the distribution function of the polymeric fluids, respectively. The notation $\Sigma(u) = \mu(\nabla_x u + \nabla_x^\top u) + \xi \div u \cdot \text{Id}$ stands for the stress tensor, with the symbol $\text{Id}$ indicating the identity matrix, and $\mu$ and $\xi$ being the viscosity coefficients satisfying the relation $2 \mu + \xi >0$. The symbol $\otimes$ denotes the tensor product for two vectors with entries $(a \otimes b)_{ij} = a_i b_j$ for $i,\ j \in \{1,2,3\}$. We assume that the pressure obeys the so-called $\gamma\mbox{-}$law, i.e. $P(\varrho) = a \varrho^\gamma$ with $\gamma \ge 1,\ a >0$. The function $U=U(q)$ is the elastic potential with $\nabla_q U$ being the spring force, $\sigma$ is a constant related to the temperature of the polymeric fluid by the relation $\sigma = k T$, where $k$ is the Boltzmann constant and $T$ is the absolute temperature. Furthermore, $r > 0$ is related to the linear damping mechanism in dynamics of the microscopic variable $q$, and $\lambda>0$ is some parameter describing the ratio between kinetic and elastic energy.

The parameter $De$ denotes the Deborah number, which represents the ratio of the time scales for elastic stress relaxation, so it characterizes the fluidity of the system. The smaller the Deborah number is, the system behaves more like a Newtonian fluid. On the other hand, the Mach number $Ma$ describes the ratio between the fluid velocity and the sound speed, hence measuring the compressibility of the system. In fact, as $De$ and $Ma$ go to zero simultaneously at the same rate, the system \eqref{eq:C-polymeric} will converges to the inhomogeneous incompressible Navier-Stokes equations. This limit will be studied in a forthcoming paper. In \cite{JLL-2016-Deborah}, the zero Deborah number limit from the incompressible version of \eqref{eq:C-polymeric} was justified. Since both of Deborah and Mach numbers are irrelevant parameters in studying the well-posedness of the model, we can thus set their values to be 1 without loss of generality in the rest of this paper.

The first two equations in system \eqref{eq:C-polymeric} are the compressible Navier-Stokes equations, with the forcing term coupled with the third equation which is a nonlinear Fokker-Planck equation. They are coupled through the induced elastic stress term in the fluid equation reflecting the micro-macro interaction, and the macro-micro drift term in the kinetic equation, corresponding to the last term on the right-hand side of \eqref{eq:C-polymeric}$_2$ and the third term on the left-hand side of \eqref{eq:C-polymeric}$_3$, respectively. More importantly, this micro-macro coupling system obeys the variational structure, and can be derived from the basic energy law by the energetic variational method (see \S 2 below):
	\begin{align}\label{eq:basic energy}
		&\tfrac{\d}{\d t} \left\{
				\int_\Omega \left[ \tfrac{1}{2} \varrho |u|^2 + \tfrac{a \varrho^\gamma}{\text{Ma}^2(\gamma-1)} \right]  \d x
				+ \tfrac{\lambda}{\text{De}} \int_\Omega\int_{\mathbb{R}^3} \Psi \left( \sigma \ln \Psi + \tfrac{1}{r} U \right)  \d q \d x
			\right\}
	\\\no
		+ & \int_\Omega \left[ \mu |\nabla_x u|^2 + (\mu+\xi)|\div u|^2 \right] \d x
		+ \tfrac{\lambda}{\text{De}^2} \int_\Omega \int_{\mathbb{R}^3} \Psi \left|\nabla_q \left( \sigma \ln \Psi + \tfrac{1}{r} U \right)  \right|^2 \d q \d x
	= 0.
	\end{align}

\subsection{Brief review of related research} % (fold)
\label{sub:brief_review_of_related_research}

% subsection brief_review_of_related_research (end)

There are two important models in the studying of viscoelastic fluids: one is the FENE (Finite Extensible Nonlinear Elastic) model and the other is Oldroyd-B model, up to the difference of potential function. In the simplest case of the Hookean spring, $U(q) = \tfrac{1}{2} |q|^2$ with $q$ belonging to the whole space $\mathbb{R}^3$, which leads to a closed equation for the extra stress tensor term $\int_{\mathbb{R}^3} q \otimes \nabla_q U \Psi \d q$ (second moment of $\Psi$), then we can recover the Oldroyd-B model. While for the FENE model, the polymer elongation vector is commonly assumed to be in a bounded open ball $B(0, b_0)$ of radius $b_0 >0$, where the potential $U(q) = - k \ln (1- \tfrac{|q|^2}{|b_0|^2}) $ for some constant $k>0$. Both of two models have been extensively studied, and most of researches are concentrated on the (approximate) closure procedure of the Fokker-Planck equation. The local existence of those systems were proved by many researchers in different setting, see for instance \cite{ELZ-2004-CMP,JLL-2004-JFA,Masmdi-2008-CPAM,ZZ-2006-ARMA}. Chemin-Masmoudi \cite{CM-2001-SIMA} proved in critical Besov spaces the local and global well-posedness. We also refer the reader to Lions-Masmoudi \cite{LM-2000-CAM} concerning the global weak solutions. Lei-Liu-Zhou \cite{LLZ-2008-ARMA} proved the global existence for small data. Lin-Zhang-Zhang \cite{LZZ-2008-CMP} and Masmoudi \cite{Masmdi-2008-CPAM} proved the same global results for FENE model. Moreover, Masmoudi \cite{Masmdi-2013-Invent} proved the global weak solutions to the FENE model by finding some new a priori estimates. Lin-Liu-Zhang \cite{LLZ-2007-CPAM} studied a incompressible micro-macro polymeric system and proved global existence near equilibrium with some assumptions on the potential $U$ (which contains Oldroyd-B model case). In what follows, we will also make these assumptions on $U$ for simplicity. On the other hand, there are many researches on other different micro-macro (or, kinetic-fluid) coupled models, which involve incompressible and compressible fluids, and investigate the existence of weak or strong solutions, large time behaviors and decay rates, stability, asymptotic analysis and related problems. We refer the readers to such as \cite{CG-2006-CPDE,CFTZ-2007-CMP,CM-2008-CMP,MV-2006-CPDE,OT-2008-CMP} and references therein.

In this present paper, we consider the case that polymer particles are transported by a compressible velocity fluid field, due to the fact that a compressible fluid flow is more realistic and complete, compared to its corresponding incompressible fluid flow model. In fact, the exact formulation of the micro-macro model \eqref{eq:C-polymeric} for the compressible polymeric fluid can be derived by the energetic variational approach, based on the basic energy law \eqref{eq:basic energy} obeyed by the fluid. Specifically speaking, the compressible model \eqref{eq:C-polymeric} possesses a complete conservative formulation. Furthermore, the system \eqref{eq:C-polymeric} contains more information about scales. We also emphasize that the scales represented by the Mach number and the Deborah number will play important roles in asymptotic analysis for the compressible polymeric model \eqref{eq:C-polymeric}, as we discussed before. Our aim here is to study the global existence of classical solutions to the compressible polymeric model \eqref{eq:C-polymeric} near equilibriums with small initial data. Moreover, the pressure $P(\varrho) = a \varrho^\gamma$ is assumed to be smooth enough, and the coefficients are assumed to one for simplicity.

For the microscopic equation \eqref{eq:C-polymeric}$_3$, it follows that the global equilibrium state can be defined by the Maxwellian $M(q) = \tfrac{e^{-U}}{\int_{\mathbb{R}^3} e^{-U}\d q}$, and furthermore, we can assume that $\int_{\mathbb{R}^3} e^{-U}\d q =1$ after renormalization. Then we consider the solution $(\varrho,u,\Psi)$ near the global equilibrium $(1,\, 0,\, M)$, i.e.
	\begin{align}\label{eq:perturbation}
	  \varrho = 1+ \rho, \quad u =u, \quad \Psi= M(1 + g).
	\end{align}
Inserting the above expansion into the model \eqref{eq:C-polymeric} yields us to get one new system for perturbations $(\rho,\, u,\, g)$:
	\begin{align}\label{eq:C-polymeric-pert}
	  \begin{cases}
	  	\p_t \rho + u \nabla_x \rho + (1+\rho) \div u = 0,
		\\[5pt]
			\p_t u + u \cdot \nabla_x u + \tfrac{P'(1+\rho)}{1+\rho} \nabla_x \rho
				= \tfrac{1}{1+\rho} \div \Sigma(u) + \tfrac{1}{1+\rho} \div \int_{\mathbb{R}^3} \nabla_q U \otimes q g M \d q ,
	  \\[7pt]
	  	\p_t g + u\cdot \nabla_x g + \nabla_x u q \nabla_q g + \mathcal{L} g
		 	= - 2(1+g) \div u + (1+g) \nabla_x u q \nabla_q U ,
	  \end{cases}
	\end{align}
where $\L g = - \tfrac{1}{M} \nabla_q \cdot ( \nabla_q \Psi + \nabla_q U \Psi) $ is the linear operator. Indeed, we can infer from performing the above expansion that $\L g = - \tfrac{1}{M} \nabla_q \cdot ( M \nabla_q g)$.

Before representing our main result in this paper, we first introduce some notations.

\bigskip
\noindent{\bf Notations.} For notational simplicity, we denote by $\skp{\cdot}{\cdot}$ the usual $L^2\mbox{-}$inner product in variables $x$, by $|\cdot|_{L^2_x}$ its corresponding $L^2$-norm, and by $|\cdot|_\hx{s}$ the higher order derivatives $H^s\mbox{-}$norm in variables $x$.

On the other hand, when we consider the spatial variables $x$ and microscopic variables $q$ at the same time, it is convenient for us to introduce the weighted Sobolev spaces. For that, we denote by $\langle \cdot,\, \cdot \rangle$ the standard $L^2$ inner product in both variables $x$ and $q$, and by $\nm{\cdot}_{L^2_{x,q}}$ its corresponding norm, then we can define the weighted $L^2$ inner product that
	\begin{align*}
		\skpm{f}{g} \triangleq \left\langle f,\, g M \right\rangle
			= \iint_{\Omega \times \mathbb{R}^3} fg \,M \d q\d x,
	\end{align*}
for any pairs $f(x,\, q), g(x,\, q) \in L^2_{x,\, q}$. We also use $\nm{\cdot}_\lxq$ to denote the weighted $L^2$-norm with respect to the measure $M \d q\d x$.

Let $\alpha=(\alpha_1,\, \alpha_2,\, \alpha_3) \in \mathbb{N}^3 $ be a multi-index with its length defined as $\textstyle |\alpha| = \sum_{i=1}^3 \alpha_i $. We also define the multi-derivative operator $\nabla^\alpha_x = \p_{x_1}^{\alpha_1} \p_{x_2}^{\alpha_2} \p_{x_3}^{\alpha_3}$, sometimes we also denote $\nabla^k$ for $|\alpha| = k$. In the following texts, we will use frequently the notation $\nabla^\alpha_\beta \triangleq \nabla^\alpha_x \nabla^\beta_q$ to stand for the mixed derivatives over the variables $x$ and $q$, where $\alpha=(a_1,\, a_2,\, a_3),\, \beta=(b_1,\, b_2,\, b_3)$ are multi-indices in $\mathbb{R}^3$ satisfying $|\alpha|= |a_1|+|a_2|+|a_3|$ and $|\beta|= |b_1|+|b_2|+|b_3|$.

Now we introduce the following mixed higher order derivatives,
	\begin{align*}
	  \nm{f}^2_\hxq{s} =\ & \sum_{|\alpha|\le s} \iint_{\Omega \times \mathbb{R}^3} |\nabla^\alpha_x f|^2 M \,\d q\d x,
	  \\
	  \nm{f}^2_\hxqm{s} =\ & \sum_{|\alpha|+|\beta|\le s} \iint_{\Omega \times \mathbb{R}^3} |\nabla^\alpha_\beta f|^2 M \,\d q\d x.
	\end{align*}
Note that we have $\hxqm{0} = \lxq$.

Let $\agl{q} = (1+|q|^2)^{1/2}$, then we introduce the energy and energy dissipation functionals for the fluctuation $(\rho, u, g)$ as follows:
	\begin{align}
	  E(t)=\ & |\rho|^2_\hx{3} + |u|^2_\hx{3} + \nm{\agl{q} g}^2_\hxqm{3},
	\\
	  D(t)=\ & ( \mu|\nabla_x u|^2_\hx{3} + (\mu+\xi)|\div u|^2_\hx{3} )
	  	+ \nm{\agl{q} \nabla_q g}^2_\hxqm{3}.
	\end{align}

At last, we mention that the notation $A \ls B$ will be used in the following texts to indicate that there exists some constant $C>0$ such that $A \le C B$. Furthermore, the notation $A \sim B$ means that the terms of both sides are equivalent up to a constant, namely, there exists some constant $C$ such that $ C^{-1} B \le A \le C B$.

%% ------------------------------------------------------------------------ %%

\subsection{Main result} % (fold)
\label{sub:main_result}

% subsection main_result (end)

Our aim is to establish the global-in-time existence of classical solutions to the micro-macro compressible polymeric systems \eqref{eq:C-polymeric} if the initial data are close to the global equilibrium $(1,0,M)$ in some appropriate spaces. To avoid some complicated technical treatment, we make here the same assumptions on the potential $U$ as that of \cite{LLZ-2007-CPAM}, more precisely, assume that
	\begin{align}\label{asmp-1}
	  & |q| \ls (1+ |\nabla_q U|), \quad \text{(sometimes we just assume $|q| \ls |\nabla_q U|$ for simplicity)}, \no\\
	  & \Delta_q U \le C + \delta |\nabla_q U|^2 \quad \text{ for } \delta < 1, \\\no
	  \int_{\mathbb{R}^3} & |\nabla_q U|^2 M \d q \le C, \quad \int_{\mathbb{R}^3} |q|^4  M \d q \le C,
	\end{align}
and
	\begin{align}\label{asmp-2}
	  &	|\nabla_q^k(q \nabla_q U)| \ls (1+|q| |\nabla_q U|), \no\\
	  &	\int_{\mathbb{R}^3} |\nabla_q^k(q \nabla_q U \sqrt{M})|^2 \d q \le C, \\\no
	  &	\abs{\nabla_q^k(\Delta_q U - \tfrac{1}{2} |\nabla_q U|^2) } \ls (1+ |\nabla_q U|^2),
	\end{align}
with the integer $1\le k \le 3$.

Now we state the main result of this paper on the global-in-time existence result of classical solutions to the compressible polymeric systems \eqref{eq:C-polymeric} near equilibrium.
\begin{theorem}[Global existence]\label{thm:main}
	Let $(\rho,u,g)$ defined through \eqref{eq:perturbation} be the fluctuation near the global equilibrium $(1,0,M)$ of the compressible polymeric system \eqref{eq:C-polymeric}, with their initial data $(\rho_0,u_0,g_0)$ satisfying the conditions $\Psi_0= M(1 + g_0)>0$ and $\int_{\mathbb{R}^3} \Psi_0 \d q =1$.

	Then, there exists some constant $\ep$ sufficiently small, such that, if the initial fluctuation satisfies 
		\begin{align}
		  & E(0) = |\rho_0|^2_\hx{3} + |u_0|^2_\hx{3} + \nm{\agl{q} g_0}^2_\hxqm{3} \le \ep, \\
		  & \int_\Omega \left[ \tfrac{1}{2} \rho_0 |u_0|^2 + \tfrac{a \rho_0^r}{(r-1)} \right]  \d x
				+ \int_\Omega\int_{\mathbb{R}^3} \Psi_0 (\ln \Psi_0 + U) \d q \d x \le \ep,
		\end{align}
	then the compressible polymeric system \eqref{eq:C-polymeric} admits a unique global classical solution $(\varrho,u,\Psi)$ with $\Psi = M(1 + g)>0$, and moreover,
		\begin{align}
		  \sup_{t \in [0,+\infty)} E(t) + \int_0^{+\infty} D(t) \d t \le \ep.
		\end{align}

\end{theorem}

To prove Theorem \ref{thm:main}, we need firstly to prove the existence of local solutions in some appropriate spaces. Introducing $\mathcal{E}(t)= E(t) + \int_0^t D(s)\d s$ with $\mathcal{E}(0) = E(0)$, we state the following local existence result, whose proof will be given via a standard iterating method (see Section 4 below):
\begin{proposition}[Local existence] \label{prop:local exis}
	Assume $E(0) \le M_0/2 $ for some constant $M_0>0$, then there exists a time $T_*>0$ such that the compressible polymeric system \eqref{eq:C-polymeric-pert} admits a unique local classical solution $(\rho,u,g) \in L^\infty(0, T_*; H^3_x \times H^3_x \times \hxqm{3}) $, moreover,
		\begin{align}
		  \E(t) \le M_0.
		\end{align}
		
\end{proposition}

The most important part of proving Theorem \ref{thm:main} is to get the a priori estimate for local solutions, which only depends on the initial datum. We present that in the following proposition.
\begin{proposition}[A priori estimate] \label{prop:A-priori}
	Let the triples $(\rho(t,x),u(t,x),g(t,x,q))$ be the local solutions constructed in Proposition \ref{prop:local exis}. Then there exist $\ep>0$ and $C_0>1$, such that if $\sup_{t \in [0, T_*]} E(t) \le \ep$, then it holds that
	\begin{align}
	  E(t) \le C_0 E(0).
	\end{align}
	
\end{proposition}

Briefly speaking, Proposition \ref{prop:A-priori} is achieved by estimating the higher order derivatives of fluctuations for density, velocity, and the microscopic probability distribution function, over the pure spatial variable $x$ and the mixed macro-micro variables $(x,q)$. Then combining the usual energy methods of Matsumura-Nishida \cite{MN-1983-CMP} for compressible Navier-Stokes equations, and noticing the cancellation relation between equations \eqref{eq:deriv u} and \eqref{eq:deriv g} (see \eqref{eq:cancellation} below), we can deduce the closed a priori estimate for the fluctuation system \eqref{eq:C-polymeric-pert}. The procedure will be separated into several subsections in Section 3, where the assumptions on potential function $U$ \eqref{asmp-1}-\eqref{asmp-2} and Lemma \ref{lemm:weight inequ} coming from the weighted Poincar\'e inequality \eqref{eq:weight Poinc} will play important roles. More details will be given there.

The plan of this paper is as follows: In the sequel we present some preliminaries which will be useful for our later proof. In Section 2 the exact formulation of micro-macro model for the compressible polymeric fluid \eqref{eq:C-polymeric} is formally derived, via the energetic variational approach. Section 3 is devoted to obtain the uniform a priori estimate stated in Proposition \ref{prop:A-priori} by closing the higher order derivatives estimates for fluctuations, including both pure macroscopic spatial derivatives and macro-micro mixed derivatives. In the last section, we construct the local solutions by an iteration scheme, based on which we can complete the proof of global existence result in Theorem \ref{thm:main}, by combining with Proposition \ref{prop:A-priori} and a standard continuum argument.

%% ------------------------------------------------------------------------ %%

\subsection{Preliminaries} % (fold)
\label{sub:preliminaries}

% subsection preliminaries (end)

For the convenience of readers, we give some lemmas which will be frequently used in the rest of the paper. Firstly, we state the Moser-type inequality formulated by Moser \cite{Moser-1966-1}, which can be found in many other references as well, for instance \cite{Majda-1984-BOOK,Taylor-2011-BOOK}. The first concerns the commutator estimates, i.e.
	\begin{lemma}[Moser-type inequality]
	\label{lemm:Moser inequ}

		For functions $f,\, g \in H^m \cap L^\infty$, $m \in \mathbb{Z}_+ \cup\, \{ 0\} $, and $|\alpha| \le m$, we have
			\begin{align}
			  \abs{\nabla^\alpha(f g) - f \nabla^\alpha g}_{L^2} \le C \left( |\nabla^m f|_{L^2} |g|_{L^\infty} + |\nabla_x f|_{L^\infty} |\nabla^{m-1} g|_{L^2} \right),
			\end{align}
		with the constant $C$ depends only on $s$.
	\end{lemma}

Combining with the chain rules, the above Moser-type inequality also yields the following lemma on smooth compositions:
	\begin{lemma}\label{lemm:composition}
		Let $F$ be a smooth function satisfying that $F(0)=0$. Then for any integer $s>0$ and $u \in H^s \cap L^\infty$, there exists constant $C$ depending only on $s$ and $|u|_{L^\infty}$, such that $F(u) \in H^s$, moreover,
			\begin{align}
			  |F(u)|_{H^s} \le C|u|_{H^s}.
			\end{align}

	\end{lemma}

By the Poincar\'e inequality, we can prove the following lemma, which is just a variation of Lemma 3.2 in \cite{LLZ-2007-CPAM} and Lemma 3.3 in \cite{JLL-2016-Deborah}. This lemma will play an important role in controlling the energy estimates for microscopic variables.
	\begin{lemma}\label{lemm:weight inequ}
		Assume $g \in \hxqm{s}$ with $\int_{\mathbb{R}^3} g M \d q= 0$, then we have
			\begin{align}\label{eq:weight inequ-U}
				\nm{\nabla_q U g}_\hxq{s} \ls\ & \nm{\nabla_q g}_\hxq{s}, \quad
			  \nm{q g}_\hxq{s} \ls \nm{\nabla_q g}_\hxq{s},
			  \\\label{eq:weight inequ-qU}
			  \nm{q \nabla_q U g}_\hxq{s} \ls\ & \nm{\agl{q} \nabla_q g}_\hxq{s}, \quad
			  \nm{|q|^2 g}_\hxq{s} \ls \nm{\agl{q} \nabla_q g}_\hxq{s},
			\end{align}
		and
			\begin{align}
			  \nm{q \nabla_q U \nabla^\alpha_\beta g}_\lxq
			    \ls\ & \nm{\agl{q} \nabla_q \nabla^\alpha_\beta g}_\lxq + \nm{\nabla^\alpha_\beta g}_\lxq , \\
			  \nm{|q|^2 \nabla^\alpha_\beta g}_\lxq
			    \ls\ & \nm{\agl{q} \nabla_q \nabla^\alpha_\beta g}_\lxq + \nm{\nabla^\alpha_\beta g}_\lxq .
			\end{align}
			
	\end{lemma}

\proof We split the proof into three steps.

\smallskip\noindent\underline{\bf Step 1.} We firstly prove the results in $\lxq$, corresponding to the case $s=0$. Noticing the assumption $\int g M \d q= 0$, the Poincar\'e inequality with weight yields immediately that
	\begin{align}\label{eq:weight Poinc}
	  \int g^2 M \d q \le C \int |\nabla_q g|^2 M \d q.
	\end{align}
	
Performing an integration by part, it follows
	\begin{align}
	  \int |\nabla_q U|^2 g^2 M \d q
	  =\ & \int \Delta_q U g^2 M \d q + 2 \int \nabla_q U g \nabla_q g M \d q \\\no
	  \le\ & \int (C+ \delta |\nabla_q U|^2) g^2 M \d q + \tfrac{1}{4} \int |\nabla_q U|^2 g^2 M \d q + 4 \int |\nabla_q g|^2 M \d q,
	\end{align}
where we have used the assumptions on $U$ \eqref{asmp-1} and the H\"older inequality. Then combining with the above weighted Poincar\'e inequality \eqref{eq:weight Poinc} and taking some small $\delta$, we get the first result in \eqref{eq:weight inequ-U} on the case $s=0$:
	\begin{align}
	  \int |\nabla_q U|^2 g^2 M \d q \le C \int |\nabla_q g|^2 M \d q,
	\end{align}
which implies immediately the second result due to the assumptions on $U$ \eqref{asmp-1}.

On the other hand, we can write
	\begin{align}
	  & \iint |q|^2 |\nabla_q U|^2 g^2 M \d q \d x \\\no
	  \ls\ & \iint |\nabla_q U|^2 \abs{q g - \int q' g(q') M \d q'}^2 M \d q \d x
	  		+ \int |\nabla_q U|^2 M \d q \cdot \int \abs{\int q' g(q') M \d q'}^2 \d x \\\no
	 	\triangleq\ & I_1 + I_2.
	\end{align}
By the above results \eqref{eq:weight inequ-U} on the case $s=0$ and \eqref{eq:weight Poinc}, we have
	\begin{align}
	  I_1 \ls \iint |\nabla_q (q g)|^2 M \d q \d x
	  		\ls \iint (1+|q|^2) |\nabla_q g|^2 M \d q \d x.
	\end{align}
By the assumptions on $U$ \eqref{asmp-1} and the H\"older inequality, we get
	\begin{align}
	  I_2 \le C \int |q|^2 M \d q \cdot \iint g^2 M \d q \d x
	  \ls \iint |\nabla_q g|^2 M \d q \d x,
	\end{align}
where we have used the weighted Poinc\'e inequality \eqref{eq:weight Poinc} again. Together with the above two inequalities, the first result \eqref{eq:weight inequ-qU} on the case $s=0$ is proved, which also yields the second inequality in \eqref{eq:weight inequ-qU}.

\smallskip\noindent\underline{\bf Step 2.} We secondly prove the results in $\hxq{s}$. Indeed, notice that the zero mean value condition $\int \nabla^k_x g M \d q= 0$ for any integer $k \le s$ remains true, we immediately get the results \eqref{eq:weight inequ-U}--\eqref{eq:weight inequ-qU}, by repeating the process in step 1 and replacing $g$ by $\nabla^k_x g$.

\smallskip\noindent\underline{\bf Step 3.} We now prove the rest results for mixed derivatives. Let $F= \nabla^\alpha_\beta g$ with $|\alpha|+|\beta|\le s$ and $|\beta|\ge 1$. A similar argument as above enables us to get
	\begin{align}
	  \int |\nabla_q U|^2 F^2 M \d q
	  \ls \int F^2 M \d q + \int |\nabla_q F|^2 M \d q,
	\end{align}
and
	\begin{align}
	  \iint |q|^2 |\nabla_q U|^2 F^2 M \d q \d x
	  \ls \iint (1+|q|^2) |\nabla_q F|^2 M \d q \d x + \iint F^2 M \d q \d x.
	\end{align}
Combining the above two inequalities leads to the penultimate result, while the last one is thus proved by the assumptions on $U$ \eqref{asmp-1}. \qed

%% ------------------------------------------------------------------------ %%

%% ------------------------------------------------------------------------ %%

\section{Derivation of compressible polymeric systems} % (fold)
\label{sec:derivation_of_compressible_polymeric_systems}

% section derivation_of_compressible_polymeric_systems (end)

In this section, we employ the energetic variational approach to formally derive the micro-macro model \eqref{eq:C-polymeric} for the compressible polymeric fluid. We first define the particle trajectory $X(\alpha,\, t)$ with respect to the Lagrangian coordinate $\alpha$, as follows,
	\begin{align}\label{eq:particle trajectory}
		\begin{cases}
		  \tfrac{\d}{\d t} X(\alpha,\, t) = u(X(\alpha,\, t),\, t),
		\\[2pt]
		  X(\alpha,\, 0) = \alpha.
		\end{cases}
	\end{align}

We then introduce the deformation tensor
	\begin{align}
	  F(X(\alpha,\, t),\, t) = \nabla_\alpha X,
	\end{align}
which satisfies the transport law
	\begin{align}
	  \partial_t F + u \cdot \nabla F = \nabla u F.
	\end{align}
In fact, the deformation tensor $F$ carries all the transport and kinematic informations of the material in complex fluid theory. To derive the exact formulation of compressible polymeric fluid system, we perform the energetic variational approach which is developed mainly by Liu and his collaborators in their series works, for instance \cite{HKL-2010-DCDS,Liu-2009-NOTES,SL-2009-DCDS}. Roughly speaking, this method mainly concerns two fundamental scalar quantities, namely the kinetic energy and the free energy, since the two scalars contain the complete dynamics of the fluid system. There are two fundamental variational principles that are usually used to model the complex system, i.e. the least action principle and the maximum dissipation principle.

\subsection{Derivation of momentum equations} % (fold)
\label{sub:derivation_of_momentum_equations}

% subsection derivation_of_momentum_equations (end)

For the momentum equation, we consider the competition between the macroscopic kinetic energy and the averaged effects due to the microscopic elastic energy. We start from the least action principle, which reads,
	\begin{align}
	  A (x) =\ & \int_0^T \int_{\Omega_t} \tfrac{1}{2} \varrho |\dot x|^2 \d x \d t
	  			  + \lambda \int_0^T \int_{\Omega_t} \int_{\mathbb{R}^3} (\sigma \Psi \ln \Psi + \tfrac{1}{r} U \Psi) \d q \d x \d t.
	\end{align}
where we mention again $\lambda$ indicates the ratio between the kinetic and elastic energy. Notice the Cauchy-Born relation (see \cite{TM-2011-BOOK} for instance) between the microscopic variable $q$ and its Lagrangian director variable $q'$ that $q=Fq'$, then it follows,
	\begin{align}
	  A (x) =\ & \int_0^T \int_{\Omega_0} \tfrac{1}{2} |\dot X|^2 \left[ \varrho(X(\alpha,\, t)) \det \nabla_\alpha	X \right] \d \alpha \d t
	  	\\\no
	  		 		&	+ \lambda \int_0^T \int_{\Omega_0} \int_{\mathbb{R}^3} (\sigma \Psi \ln \Psi + \tfrac{1}{r} U \Psi) \circ (X,q) \left( \det \nabla_\alpha	X \right)^2 \d q' \d \alpha \d t
	\\\no
				=\ & \int_0^T \!\!\int_{\Omega_0} \tfrac{1}{2} |\dot X|^2 \varrho_0(\alpha) \d \alpha \d t
	  		 		+ \lambda \int_0^T \!\!\int_{\Omega_0} \!\int_{\mathbb{R}^3} (\sigma \Psi \ln \Psi + \tfrac{1}{r} U \Psi) \circ (X,q) \left( \det \nabla_\alpha	X \right)^2 \d q' \d \alpha \d t,
	\end{align}
where we have used the compressible relation $\varrho(X(\alpha,\, t)) \det \nabla_\alpha	X = \varrho_0(\alpha)$.

We are now ready to perform the energetic variational approach. Let $X^\ep = X + \ep Y$, then we immediately get
	\begin{align}
	  \tfrac{\d}{\d \ep} X^\ep |_{\ep=0} = Y,
	  \quad \text{ and } \quad
	  \tfrac{\d}{\d \ep} \det \nabla_\alpha X^\ep |_{\ep=0} = \div Y \det \nabla_ \alpha X.
	\end{align}

Performing a variation with respect to $\ep$ reads that
	\begin{align}
	  \tfrac{\d}{\d \ep} A^\ep(x) |_{\ep=0}
	  	=\ & \int_0^T \int_{\Omega_0} \varrho_0(\alpha) \left[ \dot X^\ep \tfrac{\d}{\d \ep} \dot X^\ep |_{\ep=0} \right] \d \alpha \d t
			\\\no
			& + \lambda \int_0^T \int_{\Omega_0} \int_{\mathbb{R}^3} (\sigma (\ln \Psi + 1) + \tfrac{1}{r} U) \nabla_X \Psi \tfrac{\d}{\d \ep} X^\ep |_{\ep=0} (\det \nabla_ \alpha X)^2 \d q' \d \alpha \d t
			\\\no
			& + \lambda \int_0^T \int_{\Omega_0} \int_{\mathbb{R}^3} (\sigma (\ln \Psi + 1) + \tfrac{1}{r} U) \nabla_q \Psi \tfrac{\d}{\d \ep} F^\ep |_{\ep=0} q' (\det \nabla_ \alpha X)^2 \d q' \d \alpha \d t
			\\\no
			& + \lambda \int_0^T \int_{\Omega_0} \int_{\mathbb{R}^3} (\sigma \Psi \ln \Psi + \tfrac{1}{r} U \Psi)\circ (X,q) \ 2 \det \nabla_ \alpha X \tfrac{\d}{\d \ep} \det \nabla_\alpha X^\ep  |_{\ep=0} \d q' \d \alpha \d t
	\\\no
		 =\ & - \int_0^T \int_{\Omega_0} \varrho_0(\alpha) \dot u (X(\alpha,\, t),\, t) Y \d \alpha \d t
		 	\\\no
		 	  & + \lambda \int_0^T \int_{\Omega_0} \int_{\mathbb{R}^3} (\sigma (\ln \Psi + 1) + \tfrac{1}{r} U) \nabla_X \Psi Y (\det \nabla_ \alpha X)^2 \d q' \d \alpha \d t
			\\\no
		 	  & + \lambda \int_0^T \int_{\Omega_0} \int_{\mathbb{R}^3} (\sigma (\ln \Psi + 1) + \tfrac{1}{r} U) \nabla_q \Psi \nabla_X Y Fq' (\det \nabla_ \alpha X)^2 \d q' \d \alpha \d t
		 	\\\no
				& + 2 \lambda \int_0^T \int_{\Omega_0} \int_{\mathbb{R}^3} (\sigma \Psi \ln \Psi + \tfrac{1}{r} U \Psi)\circ (X,q) \ \div Y (\det \nabla_ \alpha X)^2 \d q' \d \alpha \d t,
	\end{align}
where we have used $\nabla_q \cdot (\nabla_x Y q) = \div Y $ and
	\begin{align*}
	  \tfrac{\d}{\d \ep} F^\ep |_{\ep=0} = \tfrac{\d}{\d \ep} \nabla_\alpha X^\ep |_{\ep=0}
	  = \nabla_ \alpha Y = \nabla_X Y \nabla_\alpha X = \nabla_X Y F.
	\end{align*}

Noticing the fundamental facts that
	\begin{align*}
	  \nabla_X (\sigma \Psi \ln \Psi + \tfrac{1}{r} U \Psi) =& [\sigma (\ln \Psi + 1) + \tfrac{1}{r} U] \nabla_X \Psi, \\\no
	  \nabla_q (\sigma \Psi \ln \Psi + \tfrac{1}{r} U \Psi) =& [\sigma (\ln \Psi + 1) + \tfrac{1}{r} U] \nabla_q \Psi + \tfrac{1}{r} \nabla_q U \Psi,
	\end{align*}
and performing an integration by part, it follows that,
	\begin{align}
		 \tfrac{\d}{\d \ep} A^\ep(x) |_{\ep=0}
		=\ & - \int_0^T \int_{\Omega_0} \left[ \varrho (\partial_t u + u \cdot \nabla u) \right] \circ X \ Y \det \nabla_\alpha X \d \alpha \d t
		 	\\\no
				& - \tfrac{\lambda}{r}  \int_0^T \int_{\Omega_0} \int_{\mathbb{R}^3} \nabla_q U \Psi \nabla_X Y q \det \nabla_\alpha X \d q \d \alpha \d t
	\\\no
		=\ & - \int_0^T \int_{\Omega_t} \left[ \varrho (\partial_t u + u \cdot \nabla u) \right] Y \d x \d t
				+ \tfrac{\lambda}{r} \int_0^T \int_{\Omega_t} \int_{\mathbb{R}^3} \div (\nabla_q U \otimes q \Psi) Y \d q \d x \d t.
	\end{align}

By postulating the viscosity and pressure contributions included in the stress tensor, we get finally
	\begin{align}
	  \varrho (\partial_t u + u \cdot \nabla_x u) + \nabla_x p
	  = \div \left[ \mu(\nabla_x u + \nabla_x^\top u) + \xi (\div u) \text{Id} \right]
	  	+ \tfrac{\lambda}{r} \div \int_{\mathbb{R}^3} (\nabla_q U \otimes q \Psi) \d q ,
	\end{align}
with $p(\varrho) = a \varrho^\gamma$. This equation will lead us to the conservative form \eqref{eq:C-polymeric}$_2$, under the condition that both of Deborah number and Mach number are equal to 1.

\subsection{Derivation of microscopic equation} % (fold)
\label{sub:derivation_of_microscopic_equation}

% subsection derivation_of_microscopic_equation (end)

Recall the stochastic model describing the dynamics of the spring $q$ that
	\begin{align}
	  \d q = - \tfrac{1}{r} \nabla_q U \d t + \sigma \d W_t,
	\end{align}
where $U(q)$ is a given elastic potential and $W_t$ is the regular Weiner process. It follows from the Ito's integration lemma that the spatially homogeneous distribution function $\Psi = \Psi(q,\, t)$ solve the following Fokker-Planck equation,
	\begin{align}
	  \partial_t \Psi = \sigma \Delta_q \Psi + \tfrac{1}{r} \nabla_q \cdot (\nabla_q U \Psi).
	\end{align}

Now we are in a position to consider the transport of distribution function $\Psi = \Psi(x,\, q,\, t)$ in spatial and spring domains $(x,\, q) \in \Omega \times \mathbb{R}^3$ in a macro-micro scales coupling system. The conversation law of total mass yields the transport formula that
	\begin{align}
	  \tfrac{\d \Psi}{\d t} = \partial_t \Psi + \div (u \Psi) + \nabla_q \cdot (\nabla_x u q \Psi),
	\end{align}
which enables us to get the spatially inhomogeneous kinetic equation that
	\begin{align}
	  \partial_t \Psi + \div (u \Psi) + \nabla_q \cdot (\nabla_x u q \Psi)
	 = \sigma \Delta_q \Psi + \tfrac{1}{r} \nabla_q \cdot (\nabla_q U \Psi).
	\end{align}

%% ------------------------------------------------------------------------ %%

\section{ A priori estimates} % (fold)
\label{sec:a_priori_estimates}

% section a_priori_estimates (end)

This section is devoted to the proof of the a priori estimates for the fluctuational compressible polymeric fluid \eqref{eq:C-polymeric-pert} presented in Proposition \ref{prop:A-priori}. The {\em a priori} estimates is achieved by estimating mainly the higher order derivatives of fluctuation functions, and the process is divided into two principal parts: contributions from the macroscopic fluid variables $(\rho, u)$ and that from the microscopic distribution function $g$. Estimating on fluctuations of density and velocity is similar as the classical energy methods of Matsumura-Nishida \cite{MN-1983-CMP}for compressible Navier-Stokes equations, where we only need to consider the pure spatial derivatives. As for the microscopic equations, we first consider the pure spatial derivatives, in which the cancellation relation (\eqref{eq:cancellation} below) will be used to control some bad term by combining equations \eqref{eq:deriv u} and \eqref{eq:deriv g} together. However, the estimates are still not closed due to some additional terms involving higher order moment. For that, we will employ the mixed derivatives estimates with respect to the macro-micro variables $(x,q)$ in some certain weighted Sobolev spaces. With the aid of some small parameter $\eta$, we finally deduce the closed a priori estimate for the fluctuational system \eqref{eq:C-polymeric-pert} for the equivalent energy functionals $E_\eta(t)$ (hence for the original energy functionals $E(t)$). In the whole proof, the assumptions on potential function $U$ \eqref{asmp-1}-\eqref{asmp-2} and Lemma \ref{lemm:weight inequ} coming from the weighted Poincar\'e inequality \eqref{eq:weight Poinc} will play important roles, and will be frequently used. The details will be given in the sequel.

Before starting the estimates, we note that the local existence of classical solutions to the fluctuational system \eqref{eq:C-polymeric-pert} is supposed to be valid, whose proof is postponed to the next section.

\subsection{Contributions from fluid variables} % (fold)
\label{sub:contributions_from_fluid_variables}

% subsection contributions_from_fluid_variables (end)

In this subsection we consider the higher order derivatives estimates on density $\rho$ and velocity $u$. Before that, we first point out the energy functionals $\abs{\rho}^2_\hx{3} + \abs{u}^2_\hx{3}$ and $\sum_{|\alpha| \le 3} \int_{\mathbb{R}^3} \left[ \tfrac{P'(1+\rho)}{1+\rho} |\nabla^\alpha_x \rho|^2 + (1+\rho) |\nabla^\alpha_x u|^2 \right] \d x$ are in fact equivalent as we study the classical solutions in a fluctuation framework. Indeed, we have the fundamental fact $|\rho|_{L^\infty_x} \ls |\rho|_\hx{2}$ by the Sobolev embedding inequality, which implies that $\rho$ is bounded from above and below as the the small condition $|\rho|_\hx{3}< \ep$ had been assumed.

\subsubsection{Estimates on density} % (fold)
\label{ssub:estimates_on_density}

% subsubsection estimates_on_density (end)

Applying the multi-derivative operator $\nabla^\alpha_x$ with $|\alpha|=3$ to equation \eqref{eq:C-polymeric-pert}$_1$, we can infer that
	\begin{align}\label{eq:deriv rho}
	  \p_t \nabla^\alpha_x \rho + u \cdot \nabla_x \nabla^\alpha_x \rho + (1+\rho) \div \nabla^\alpha_x u
	  	+ [\nabla^\alpha_x, u \cdot \nabla_x] \rho + [\nabla^\alpha_x, (1+\rho) \div] u =0.
	\end{align}
Taking $L^2_x$ inner product with $\tfrac{P'(1+\rho)}{1+\rho} \nabla^\alpha_x \rho$ enables us to derive the estimates:
	\begin{align}
	  \skp{\p_t \nabla^\alpha_x \rho}{\tfrac{P'(1+\rho)}{1+\rho} \nabla^\alpha_x \rho}
	  	=\ & \tfrac{1}{2} \tfrac{\d}{\d t}\int \tfrac{P'(1+\rho)}{1+\rho} |\nabla^\alpha_x \rho|^2 \d x
	    	- \tfrac{1}{2} \int \p_t (\tfrac{P'(1+\rho)}{1+\rho}) |\nabla^\alpha_x \rho|^2 \d x,
	  \\
    \skp{u \cdot \nabla_x \nabla^\alpha_x \rho}{\tfrac{P'(1+\rho)}{1+\rho} \nabla^\alpha_x \rho}
    	=\ & \tfrac{1}{2} \!\! \int \!\! \tfrac{P'(1+\rho)}{1+\rho} u \cdot \nabla_x |\nabla^\alpha_x \rho|^2 \d x
    	= -\tfrac{1}{2} \! \int \div (\tfrac{P'(1+\rho)}{1+\rho} u) |\nabla^\alpha_x \rho|^2 \d x\\
	  		\ls\ & |\rho|_\hx{3} |u|_\hx{3} |\nabla^\alpha_x \rho|^2_\lx,
	\end{align}
where we have used in the second inequality the simple fact
	\begin{align*}
	  |\div (\tfrac{P'(1+\rho)}{1+\rho} u)|_{L^\infty_x}
	  =\ & |\nabla_x(\tfrac{P'(1+\rho)}{1+\rho}) u|_\hx{2}
	  	+ | \tfrac{P'(1+\rho)}{1+\rho} \div u|_{L^\infty_x} \\
	  \le\ & C |\rho|_\hx{3} |u|_\hx{2} + C |\rho|_\hx{3} |\div u|_{L^\infty_x} \\
	  \le\ & C |\rho|_\hx{3} |u|_\hx{3}.
	\end{align*}

We now deal with the two commutator terms in the above equations. By the Moser-type inequality in Lemma \ref{lemm:Moser inequ}, we have
	\begin{align*}
	  |[\nabla^\alpha_x, u \cdot \nabla_x] \rho |_\lx
	  \ls\ & |u|_\hx{3} |\nabla_x \rho|_{L^\infty_x} + |\nabla_x u|_{L^\infty_x} |\nabla_x \rho|_\hx{2} \\
	  \ls\ & |\rho|_\hx{3} |u|_\hx{3},
	\end{align*}
which implies that
  \begin{align}
    \skp{[\nabla^\alpha_x, u \cdot \nabla_x] \rho}{\tfrac{P'(1+\rho)}{1+\rho} \nabla^\alpha_x \rho}
	  \ls |\rho|_\hx{3} |u|_\hx{3} |\nabla^\alpha_x \rho|_\lx.
  \end{align}
For the last term, we get similarly,
	\begin{align}
    \skp{[\nabla^\alpha_x, (1+\rho) \div] u}{\tfrac{P'(1+\rho)}{1+\rho} \nabla^\alpha_x \rho}
	  \ls |\rho|_\hx{3} |u|_\hx{3} |\nabla^\alpha_x \rho|_\lx.
	\end{align}

By noticing equation \eqref{eq:C-polymeric-pert}$_1$ yields that
	\begin{align*}
	  |\p_t \rho |_{L^\infty_x} = |u \nabla_x \rho + (1+\rho) \div u |_{L^\infty_x}
	  \ls |\rho|_\hx{3} |u|_\hx{3} + |u|_\hx{3} ,
	\end{align*}
combining with the H\"older inequality, this implies that
	\begin{align*}
	  \int \p_t (\tfrac{P'(1+\rho)}{1+\rho}) |\nabla^\alpha_x \rho|^2 \d x
	  \ls (|\rho|_\hx{3} |u|_\hx{3} + |u|_\hx{3}) |\nabla^\alpha_x \rho|^2_\lx .
	\end{align*}

Note that the lower order derivatives are similar and much easier, and some straightforward calculations enable us to derive for the case $|\alpha|=0$ that,
	\begin{align}
	  \tfrac{1}{2} \tfrac{\d}{\d t} |\rho|^2_\lx + \skp{(1+\rho) \div u}{\tfrac{P'(1+\rho)}{1+\rho} \rho}
	  \ls\ & |u|_\hx{2} |\nabla_x \rho|_\lx |\rho|_\lx.
 	\end{align}

Summing up all the estimates for $|\alpha|\le 3$, we can derive the estimates for density $\rho$ that
	\begin{align}\label{esm:deriv rho}
	   \tfrac{1}{2} \tfrac{\d}{\d t} |\rho|^2_\hx{3}
	  	+ \sum_{|\alpha|\le 3} \skp{ (1+\rho) \div \nabla^\alpha_x u}{\tfrac{P'(1+\rho)}{1+\rho} \nabla^\alpha_x \rho}
	  \ls\ & (1+ |\rho|_\hx{3}) |\rho|_\hx{3} |u|_\hx{3} |\nabla_x \rho|_\hx{2}.
	\end{align}

%% ---------------------------------------------------------------------
\subsubsection{Estimates on velocity variables} % (fold)
\label{ssub:estimates_on_velocity_variables}

% subsubsection estimates_on_velocity_variables (end)

%% ---------------------------------------------------------------------

The process to estimating velocity variables is similar to that for the density variable. Firstly we write the equation of the higher order derivatives of $u$ as follows
	\begin{align}\label{eq:deriv u}
	  \p_t \nabla^\alpha_x u & + u \cdot \nabla_x \nabla^\alpha_x u + [\nabla^\alpha_x, u \cdot \nabla_x] u
	   + \tfrac{P'(1+\rho)}{1+\rho} \nabla_x \nabla^\alpha_x \rho + [\nabla^\alpha_x, \tfrac{P'(1+\rho)}{1+\rho} \nabla_x] \rho
	   \\\no
	  =\ & \tfrac{1}{1+\rho} \div \nabla^\alpha_x \Sigma(u) + [\nabla^\alpha_x, \tfrac{1}{1+\rho} \div] \Sigma(u)
	  	\\\no
	    & + \tfrac{1}{1+\rho} \div \nabla^\alpha_x \int \nabla_q U \otimes q g M \d q
	    % \\\no
	   	  + [\nabla^\alpha_x, \tfrac{1}{1+\rho} \div] \int \nabla_q U \otimes q g M \d q.
	\end{align}
	
Next we take $L^2_x$ inner product with the quantity $(1+\rho) \nabla^\alpha_x u$, and get term by term that
	\begin{align*}
	  \skp{\p_t \nabla^\alpha_x u}{(1+\rho) \nabla^\alpha_x u}
	  	=\ & \tfrac{1}{2} \tfrac{\d}{\d t}\int (1+\rho) |\nabla^\alpha_x u|^2 \d x
	  		- \tfrac{1}{2} \int \p_t \rho |\nabla^\alpha_x u|^2 \d x,
	  \\
	  \skp{u \cdot \nabla_x \nabla^\alpha_x u}{(1+\rho) \nabla^\alpha_x u}
	  	=\ & - \tfrac{1}{2} \int \div ((1+\rho) u) |\nabla^\alpha_x u|^2 \d x \\
	  	\ls\ & (|\rho|_\hx{3} |u|_\hx{3} + |u|_\hx{3}) |\nabla^\alpha_x u|^2_\lx,
	  \\
	  \skp{[\nabla^\alpha_x, u \cdot \nabla_x] u}{(1+\rho) \nabla^\alpha_x u}
	  	\ls\ & \left( |u|_\hx{3} |\nabla_x u|_{L^\infty_x} + |\nabla_x u|_{L^\infty_x} |\nabla_x u|_\hx{2} \right) |(1+\rho) \nabla^\alpha_x u|_\lx
	  	\\
	  	\ls\ & |u|^2_\hx{3} |\nabla^\alpha_x u|_\lx,
	  \\
	  \skp{[\nabla^\alpha_x, \tfrac{P'(1+\rho)}{1+\rho} \nabla_x] \rho}{(1+\rho) \nabla^\alpha_x u}
	  	& \ls |\nabla_x \rho|^2_\hx{2} |\nabla^\alpha_x u|_\lx,
	\end{align*}
where we have used the Sobolev embedding inequality and the Moser-type inequality in Lemma \ref{lemm:Moser inequ} for the commutator estimates.

On the other hand, we turn to treat the terms on the right-hand side of equation \eqref{eq:deriv u}. The dissipative term reads that
		\begin{align}
		  \skp{\tfrac{1}{1+\rho} \div \nabla^\alpha_x \Sigma(u)}{(1+\rho) \nabla^\alpha_x u}
		  = \ & -\skp{\nabla^\alpha_x [\mu (\nabla_x u + \nabla_x^\top u) + \xi (\div u) \cdot \text{Id}]}{\nabla_x \nabla^\alpha_x u}
		  \\\no
		  = \ & - \int [\mu |\nabla_x \nabla^\alpha_x u|^2 + (\mu + \xi) |\div \nabla^\alpha_x u|^2] \d x.
		\end{align}

For the second term on the right-hand side, taking inner product as before yields that
	\begin{align}
	  \skp{[\nabla^\alpha_x, \tfrac{1}{1+\rho} \div] \Sigma(u)}{(1+\rho) \nabla^\alpha_x u}
	  = \ & \sum_{\substack{|\alpha|=|\alpha_1|+|\alpha_2|\\ |\alpha_1| \ge 1}}
	  	\skp{\nabla_x^{\alpha_1}(\tfrac{1}{1+\rho}) \div \Sigma(\nabla_x^{\alpha_2} u)}{(1+\rho) \nabla^\alpha_x u}
	  \\\no
	  \triangleq \ & \sum_{\substack{|\alpha|=|\alpha_1|+|\alpha_2|\\ |\alpha_1| \ge 1}} A_{|\alpha_1|,|\alpha_2|}.
	\end{align}
We split the sum $|\alpha_1|+|\alpha_2|= |\alpha|=3$ with $|\alpha_1| \ge 1$ into the following three cases:
	\begin{itemize}

	 	\item[-] For the case of $|\alpha_1|=3,\ |\alpha_2|=0$, we get
	 		\begin{align}
	 		  A_{3,0} = \ & \skp{\nabla^3(\tfrac{1}{1+\rho}) \div \Sigma(u)}{(1+\rho) \nabla^3 u}
	 		  				\\\no
	 		  				\le \ & \abs{\nabla^3(\tfrac{1}{1+\rho})}_\lx \abs{\div \Sigma(u)}_{L^\infty_x} \abs{(1+\rho) \nabla^3 u}_\lx
	 		  				\\\no
	 		  				\ls \ & |\nabla_x \rho|_\hx{2} |\nabla_x u|_\hx{3} |u|_\hx{3},
	 		\end{align}
	 	where we have used the Sobolev embedding inequality and Lemma \ref{lemm:composition} for the smooth function $F(\rho) \triangleq \tfrac{1}{1+\rho} - 1$ satisfying that $F(0)=0$;

	 	\item[-] For case $|\alpha_1|=2,\ |\alpha_2|=1$, we get similarly,
	 		\begin{align}
	 		  A_{2,1} = \ & \skp{\nabla^2(\tfrac{1}{1+\rho}) \div \Sigma(D u)}{(1+\rho) \nabla^3 u}
	 		  				\\\no
	 		  				\le \ & \abs{\nabla^2(\tfrac{1}{1+\rho})}_{L^4_x} \abs{\div \Sigma(D u)}_{L^4_x} \abs{(1+\rho) \nabla^3 u}_\lx
	 		  				\\\no
	 		  				\ls \ & |\nabla_x \rho|_\hx{2} |\nabla_x u|_\hx{3} |u|_\hx{3},
	 		\end{align}
	 	where we have used the Sobolev embedding inequality $|h|_{L^4_x} \ls |h|_\hx{1}$;

	 	\item[-] For case $|\alpha_1|=1,\ |\alpha_2|=2$, we get
	 		\begin{align}
	 		  A_{1,2} = \ & \skp{D (\tfrac{1}{1+\rho}) \div \Sigma(\nabla^2 u)}{(1+\rho) \nabla^3 u}
	 		  				\\\no
	 		  				\le \ & \abs{D(\tfrac{1}{1+\rho})}_{L^\infty_x} \abs{\div \Sigma(\nabla^2 u)}_{L^2_x} \abs{(1+\rho) \nabla^3 u}_\lx
	 		  				\\\no
	 		  				\ls \ & |\nabla_x \rho|_\hx{2} |\nabla_x u|_\hx{3} |u|_\hx{3},
	 		\end{align}
	 	where we have used the Sobolev embedding inequality again.   	 	
	 \end{itemize}
Noticing that the lower order derivative estimates ($|\alpha|<3$) are much easier, we can infer that, for all summation satisfying $|\alpha|\le 3$,
	\begin{align}\label{eq:commu disspa}
		\sum_{|\alpha|\le 3}
			\skp{[\nabla^\alpha_x, \tfrac{1}{1+\rho} \div] \Sigma(u)}{(1+\rho) \nabla^\alpha_x u}
		\ls \ & |\nabla_x \rho|_\hx{2} |u|_\hx{3} |\nabla_x u|_\hx{3}.
	\end{align}
	
As the end of this subsection, we consider the last commutator term in the right-hand side of equation \eqref{eq:deriv u}. We denote that
	\begin{align}
	  I = \ & \abs{[\nabla^\alpha_x, \tfrac{1}{1+\rho} \div] \int \nabla_q U \otimes q g M \d q}_\lx \\\no
	  	= \ & \sum_{\substack{|\alpha|=|\alpha_1|+|\alpha_2|\\ |\alpha_1| \ge 1}}
	  		\abs{\nabla_x^{\alpha_1} (\tfrac{1}{1+\rho}) \int q \nabla_q \nabla_x \nabla_x^{\alpha_2} g M \d q}_\lx \\\no
	  	\le \ & \sum_{\substack{|\alpha|=|\alpha_1|+|\alpha_2|\\ |\alpha_1| \ge 1}}
	  		\abs{ \abs{\nabla_x^{\alpha_1} (\tfrac{1}{1+\rho})} \left( \int q^2 M \d q \right)^{1/2} |\nabla_q \nabla_x^{\alpha_2 +1} g|_{L^2_{M\!q}}}_\lx \\\no
	  	\le \ & \sum_{\substack{|\alpha|=|\alpha_1|+|\alpha_2|\\ |\alpha_1| \ge 1}}
	  		\abs{ |{\nabla_x^{\alpha_1} (\tfrac{1}{1+\rho})}| \cdot |\nabla_q \nabla_x^{\alpha_2 +1} g|_{L^2_{M\!q}}}_\lx \\\no	  		
	  	= \ & \sum_{\substack{|\alpha|=|\alpha_1|+|\alpha_2|\\ |\alpha_1| \ge 1}} I_{|\alpha_1|,|\alpha_2|},
	\end{align}
where we have used the simple fact by integrating by parts,
	\begin{align*}
	  \div \int \nabla_q U \otimes q g M \d q
	  =\ & \int \p_{q_i} U q_j \p_{x_i} g M \d q
	  = - \int q_j \p_{x_i} g \p_{q_j} M \d q \\\no
	  =\ & \int \p_{q_j}(q_j \p_{x_i} g) M \d q
	  = \int q \nabla_q \nabla_x g M \d q,
	\end{align*}
since $\int \delta_{ij} \p_{x_i} g M \d q = \nabla_x \int g M \d q=0$.

By a discuss process with respect to the values of $\alpha_1$ and $\alpha_2$ similar as before, we can get the estimate
	\begin{align}
	  I \ls |\nabla_x \rho|_\hx{2} \nm{\nabla_q g}_\hxq{3},
	\end{align}
and the commutator in the last term of equation \eqref{eq:deriv u} can be bounded as follows,
	\begin{align}\label{eq:commu micro}
		\skp{[\nabla^\alpha_x, \tfrac{1}{1+\rho} \div] \int \nabla_q U \otimes q g M \d q}{(1+\rho) \nabla^\alpha_x u}
			\ls |\nabla_x \rho|_\hx{2} |u|_\hx{3} \nm{\nabla_q g}_\hxq{3}.
	\end{align}

Finally, combining all these above estimates enables us to get the higher order derivatives estimates for equation \eqref{eq:deriv u} with respect to the summation $|\alpha|\le 3$, that is,
	\begin{align}\label{esm:deriv u}
	  & \tfrac{1}{2} \tfrac{\d}{\d t} |u|^2_\hx{3} + \left( \mu|\nabla_x u|^2_\hx{3} + (\mu+\xi)|\div u|^2_\hx{3} \right)
	  		+ \sum_{|\alpha|\le 3} \skp{\tfrac{P'(1+\rho)}{1+\rho} \nabla_x \nabla^\alpha_x \rho}{(1+\rho) \nabla^\alpha_x u}
	 \\\no
	  \ls \ & (1+|\rho|_\hx{3}) |u|^3_\hx{3} + |\nabla_x \rho|^2_\hx{2} |u|_\hx{3}
	  			  + |u|_\hx{3} |\nabla_x \rho|_\hx{2} |\nabla_x u|_\hx{3}
	  			  + |u|_\hx{3} |\nabla_x \rho|_\hx{2} \nm{\nabla_q g}_\hxq{3} \\\no
	  			& + \sum_{|\alpha|\le 3} \skp{\div \nabla^\alpha_x \int \nabla_q U \otimes q g M \d q}{\nabla^\alpha_x u}.
	\end{align}

%% ---------------------------------------------------------------------
\subsubsection{Estimates for fluid variables} % (fold)
\label{ssub:estimates_for_fluid_variables}

% subsubsection estimates_for_fluid_variables (end)

%% ---------------------------------------------------------------------

Note that the following cancellation fact:
	\begin{align}
	  & \skp{ (1+\rho) \div \nabla^\alpha_x u}{\tfrac{P'(1+\rho)}{1+\rho} \nabla^\alpha_x \rho}
	  		+ \skp{\tfrac{P'(1+\rho)}{1+\rho} \nabla_x \nabla^\alpha_x \rho}{(1+\rho) \nabla^\alpha_x u} \\\no
	  =\ & \skp{\div \nabla^\alpha_x u}{P'(1+\rho) \nabla^\alpha_x \rho}
	  	+ \skp{P'(1+\rho) \nabla_x \nabla^\alpha_x \rho}{\nabla^\alpha_x u} \\\no
	  =\ & - \skp{P''(1+\rho) \nabla_x \rho \nabla^\alpha_x \rho}{\nabla^\alpha_x u} \\\no
	  \ls\ & |\nabla_x \rho|_{L^\infty_x} |\nabla^\alpha_x \rho|_\lx |\nabla^\alpha_x u|_\lx,
	\end{align}
combining with the above estimates \eqref{esm:deriv rho} and \eqref{esm:deriv u}, we derive finally the higher order derivatives estimates on macroscopic fluctuations that,
	\begin{align}\label{esm:fluid}
		 & \tfrac{1}{2} \tfrac{\d}{\d t} (|\rho|^2_\hx{3} + |u|^2_\hx{3}) + ( \mu|\nabla_x u|^2_\hx{3} + (\mu+\xi)|\div u|^2_\hx{3} )
	 \\\no
	  \ls \ & (1+|\rho|_\hx{3}) (|\rho|_\hx{3} |u|_\hx{3} |\nabla_x \rho|_\hx{2}+ |u|^3_\hx{3}) + |u|_\hx{3} |\nabla_x \rho|_\hx{2} (|\nabla_x \rho|_\hx{2} + |\nabla_x u|_\hx{3} + \nm{\nabla_q g}_\hxq{3}) \\\no
	  			& + \sum_{|\alpha|\le 3} \skt{\div \nabla^\alpha_x \int \nabla_q U \otimes q g M \d q}{\nabla^\alpha_x u}.
	\end{align}

%% ---------------------------------------------------------------------
\subsubsection{Estimates on the density dissipation} % (fold)
\label{ssub:estimates_on_the_density_dissipation}

% subsubsection estimates_on_the_density_dissipation (end)

%% ---------------------------------------------------------------------

We have to deal with the contribution coming from the density dissipation. Applying the operator $\nabla^\alpha_x$ ($|\alpha| \le 2$) directly to equation \eqref{eq:C-polymeric-pert}$_2$, and taking integration  with $\nabla_x \nabla^\alpha_x \rho$, we get that,
	\begin{align}
	  & \skp{\nabla^\alpha_x \p_t u}{\nabla_x \nabla^\alpha_x \rho} + |\nabla_x \nabla^\alpha_x \rho|^2_\lx \\\no
	  =\ & \skp{\nabla^\alpha \left\{ \tfrac{1}{1+\rho} \div \Sigma(u) - u \cdot \nabla_x u - \left[ \tfrac{P'(1+\rho)}{1+\rho}-1 \right]  \nabla_x \rho \right\} }{\nabla_x \nabla^\alpha_x \rho} \\\no
	  & + \skp{\nabla^\alpha \{ \tfrac{1}{1+\rho} \div \int_{\mathbb{R}^3} \nabla_q U \otimes q g M \d q \} }{\nabla_x \nabla^\alpha_x \rho} \\\no
	  =\ & I_1 + I_2.
	\end{align}
Performing direct calculations yields that the right-hand side terms can be bounded by
	\begin{align}
		I_1 \le \frac{1}{8} |\nabla_x \nabla^\alpha_x \rho|^2_\lx + C |u|^2_\hx{2} |\nabla_x u|^2_\hx{2} + C |\rho|^2_\hx{2} |\nabla_x \rho|^2_\hx{2} + C |\rho|^2_\hx{2} |\nabla_x u|^2_\hx{2},
	\end{align}
and
	\begin{align}
		I_2 \le \frac{1}{8} |\nabla_x \nabla^\alpha_x \rho|^2_\lx
						+ C |\rho|^2_\hx{2} \nm{\nabla_q g}^2_\hxq{3}.
	\end{align}

By integration by part, we have
	\begin{align}
	  \skp{\nabla^\alpha_x \p_t u}{\nabla_x \nabla^\alpha_x \rho}
	  =\ & \tfrac{\d}{\d t} \skp{\nabla^\alpha_x u}{\nabla_x \nabla^\alpha_x \rho} - \int \div \nabla^\alpha_x u \nabla^\alpha_x \left[ u \cdot \nabla_x \rho + (1+\rho) \div u \right] \\\no
	  \ge\ & \tfrac{\d}{\d t} \skp{\nabla^\alpha_x u}{\nabla_x \nabla^\alpha_x \rho} - C (1 + |\rho|_\hx{2}) |\nabla_x u|^2_\hx{2} + C |\nabla_x \rho|_\hx{2} |\nabla_x u|^2_\hx{2},
	\end{align}
where we have used equation \eqref{eq:C-polymeric-pert}$_1$.

Thus we can infer finally that,
	\begin{align}\label{esm:density dispa}
	  & \tfrac{\d}{\d t} \Sigma_{|\alpha| \le 2} \skp{\nabla^\alpha_x u}{\nabla_x \nabla^\alpha_x \rho}
	  	+ \tfrac{4}{5} |\nabla_x \rho|^2_\hx{2} \\\no
	  \ls\ & (1 + |\rho|_\hx{2}) |\nabla_x u|^2_\hx{2} + |\nabla_x \rho|_\hx{2} |\nabla_x u|^2_\hx{2}
	   + |u|^2_\hx{2} |\nabla_x u|^2_\hx{2} + |\rho|^2_\hx{2} |\nabla_x \rho|^2_\hx{2} \\\no
	  & + |\rho|^2_\hx{2} |\nabla_x u|^2_\hx{2} + |\rho|^2_\hx{2} \nm{\nabla_q g}^2_\hxq{3}.
	\end{align}

%% ------------------------------------------------------------------------ %%
%% ------------------------------------------------------------------------ %%
\subsection{Estimates for microscopic variables} % (fold)
\label{sub:estimates_for_microscopic_variables}

% subsection estimates_for_microscopic_variables (end)

%% ---------------------------------------------------------------------

In this section we aim to estimate the microscopic equation \eqref{eq:C-polymeric-pert}$_3$, namely,
	\begin{align*}
	  \p_t g + u\cdot \nabla_x g + \nabla_x u q \nabla_q g + \mathcal{L} g
		 	= - 2(1+g) \div u + (1+g) \nabla_x u q \nabla_q U ,
	\end{align*}
where $\L g \triangleq - \tfrac{1}{M} \nabla_q \cdot ( M \nabla_q g)$ indicates the microscopic dissipation. Note that $\L$ is an adjoint operator with respect to the weighted inner product $L^2(M\!\d q)$, more precisely,
	\begin{align}
		\left\langle \L g, g \right\rangle_{L^2(M\!\d q)}
		= - \int \tfrac{1}{M} \nabla_q \!\cdot\! (M \nabla_q g) g M \!\d q
		= \int |\nabla_q g|^2 M \! \d q
		= \left\langle \nabla_q g , \nabla_q g \right\rangle_{L^2(M\!\d q)}.
	\end{align}

We will work in the weighted Sobolev space with respect to the measure $M \d q \d x$. For the sake of simplicity, we introduce one notation $\agl{q} = (1+|q|^2)^{1/2}$.

%% ---------------------------------------------------------------------
\subsubsection{Estimates for spatial derivatives} % (fold)
\label{ssub:estimates_for_spatial_derivatives}

% subsubsection estimates_for_spatial_derivatives (end)

%% ---------------------------------------------------------------------

Applying $\nabla^\alpha_x$ ($|\alpha| \le 3$) to equation \eqref{eq:C-polymeric-pert}$_3$, we get
	\begin{align}\label{eq:deriv g}
	  \p_t \nabla^\alpha_x g + \nabla^\alpha_x (u\cdot \nabla_x g) &+ \nabla^\alpha_x (\nabla_x u q \nabla_q g) + \mathcal{L} \nabla^\alpha_x g \\\no
		 	& = - \nabla^\alpha_x [2(1+g) \div u + (1+g) \nabla_x u q \nabla_q U] ,
	\end{align}
where we have used the commutative formula $\nabla^\alpha_x \mathcal{L} = \mathcal{L} \nabla^\alpha_x$.

Multiplying both sides of equation \eqref{eq:deriv g} by $\agl{q}^2 \nabla^\alpha_x g$, and integrating with respect to the weighted measure $M \d q \d x$, we firstly get
	\begin{align*}
	  \skpm{\p_t \nabla^\alpha_x g}{\agl{q}^2 \nabla^\alpha_x g}
	  =\ & \skpm{\p_t \nabla^\alpha_x g}{\nabla^\alpha_x g} + \skpm{\p_t \nabla^\alpha_x g}{|q|^2 \nabla^\alpha_x g}
	  =\tfrac{1}{2} \tfrac{\d}{\d t} \nm{\agl{q} \nabla^\alpha_x g}^2_\lxq ,
	\end{align*}
and 	
	\begin{align*}
	  \skpm{\nabla^\alpha_x (u\cdot \nabla_x g)}{\agl{q}^2 \nabla^\alpha_x g}
	  	=\ & \skpm{u\cdot \nabla_x \nabla^\alpha_x g}{\agl{q}^2 \nabla^\alpha_x g}
	  		+ \skpm{[\nabla^\alpha_x, u\cdot \nabla_x] g}{\agl{q}^2 \nabla^\alpha_x g} \no\\
	  	\ls\ & |\div u|_{L^\infty_x} \nm{\agl{q} \nabla^\alpha_x g}^2_\lxq
	  		+ |u|_\hx{3} \nm{\agl{q} \nabla_x g}_\hx{2} \nm{\agl{q} \nabla^\alpha_x g}_\lxq \no\\
	  	\ls\ & |u|_\hx{3} \nm{\agl{q} g}^2_\hxq{3},
	\end{align*}
where we have used the Sobolev inequality and the Moser-type inequality in Lemma \ref{lemm:Moser inequ} to control the commutator, as we did before.

As for the third term, we write that
	\begin{align}
	  J \triangleq \skpm{\nabla^\alpha_x (\nabla_x u q \nabla_q g)}{\agl{q}^2 \nabla^\alpha_x g}
	  =\ & \sum_{|\alpha|=|\alpha_1|+|\alpha_2|} \skpm{\nabla^{\alpha_1 +1}_x u q \nabla_q \nabla^{\alpha_2}_x g}{\agl{q}^2 \nabla^\alpha_x g} \\\no
	  \triangleq\ & \sum_{|\alpha|=|\alpha_1|+|\alpha_2|} J_{|\alpha_1|,|\alpha_2|}.
	\end{align}

Then by a similar discuss process as before, we can infer that
	\begin{align}
	  J_{3,0} \le\ & |\nabla^3_x \nabla_x u|_\lx \nm{|q| \nabla_q g}_{L^\infty_x L^2_{M\!q}} \nm{\agl{q}^2 \nabla^\alpha_x g}_\lxq \\\no
	  \ls\ & |\nabla_x u|_\hx{3} \nm{|q| \nabla_q g}_\hxq{2} \nm{\agl{q}^2 \nabla^\alpha_x g}_\lxq, \\\no
	  J_{2,1} \le\ & |\nabla^2_x \nabla_x u|_\lx \nm{|q| \nabla_q \nabla^1_x g}_{L^\infty_x L^2_{M\!q}} \nm{\agl{q}^2 \nabla^\alpha_x g}_\lxq \\\no
	  \ls\ & |u|_\hx{3} \nm{|q| \nabla_q g}_\hxq{3} \nm{\agl{q}^2 \nabla^\alpha_x g}_\lxq, \\\no
	  J_{1,2} \le\ & |\nabla^1_x \nabla_x u|_{L^4_x} \nm{|q| \nabla_q \nabla^2_x g}_{L^4_x L^2_{M\!q}} \nm{\agl{q}^2 \nabla^\alpha_x g}_\lxq \\\no
	  \ls\ & |u|_\hx{3} \nm{|q| \nabla_q g}_\hxq{3} \nm{\agl{q}^2 \nabla^\alpha_x g}_\lxq, \\\no
	  J_{0,3} \le\ & |\nabla_x u|_{L^\infty_x} \nm{|q| \nabla_q \nabla^3_x g}_\lxq \nm{\agl{q}^2 \nabla^\alpha_x g}_\lxq \\\no
	  \ls\ & |u|_\hx{3} \nm{|q| \nabla_q g}_\hxq{3} \nm{\agl{q}^2 \nabla^\alpha_x g}_\lxq,
	\end{align}
which, together with Lemma \ref{lemm:weight inequ}, implies that for $|\alpha| \le 3$,
	\begin{align}
	  J =\ & \skpm{\nabla^\alpha_x (\nabla_x u q \nabla_q g)}{\agl{q}^2 \nabla^\alpha_x g} \\\no
	  \ls\ & |\nabla_x u|_\hx{3} \nm{|q| \nabla_q g}_\hxq{2} \nm{\agl{q}^2 \nabla^\alpha_x g}_\lxq + |u|_\hx{3} \nm{|q| \nabla_q g}_\hxq{3} \nm{\agl{q}^2 \nabla^\alpha_x g}_\lxq \\\no
	  \ls\ & |\nabla_x u|_\hx{3} \nm{\agl{q} \nabla_q g}_\hxq{2} \nm{\agl{q} \nabla_q g}_\hxq{3}
	  + |u|_\hx{3} \nm{\agl{q} \nabla_q g}^2_\hxq{3}.
	\end{align}

Considering the micro-dissipation contributions term, it follows from the adjoint property of operator $\L$ that
	\begin{align}
	  \skpm{\L \nabla^\alpha_x g}{\agl{q}^2 \nabla^\alpha_x g}
	  =\ & \skpm{\L \nabla^\alpha_x g}{\nabla^\alpha_x g} + \skpm{\L \nabla^\alpha_x g}{|q|^2 \nabla^\alpha_x g} \\\no
	  =\ & \skpm{\nabla_q \nabla^\alpha_x g}{\nabla_q \nabla^\alpha_x g} - \iint \nabla_q \cdot (M \nabla_q \nabla^\alpha_x g) |q|^2 \nabla^\alpha_x g \d q\d x \\\no
	  =\ & \skpm{\nabla_q \nabla^\alpha_x g}{\nabla_q \nabla^\alpha_x g}
	  		+ \iint M \nabla_q \nabla^\alpha_x g (|q|^2 \nabla_q \nabla^\alpha_x g + 2 q \nabla^\alpha_x g) \d q\d x \\\no
	  =\ & \nm{\agl{q} \nabla_q \nabla^\alpha_x g}^2_\lxq - 3 \nm{\nabla^\alpha_x g}^2_\lxq + \iint q \nabla_q U |\nabla^\alpha_x g|^2 M\d q\d x.
 	\end{align}
Combining with the H\"{o}lder inequality and Lemma \ref{lemm:weight inequ} together yields that
	\begin{align}
	  \iint q \nabla_q U |\nabla^\alpha_x g|^2 M\d q\d x
	  \ls\ & \nm{\nabla^\alpha_x g}_\lxq \nm{q \nabla_q U \nabla^\alpha_x g}_\lxq
	  \ls \nm{\nabla^\alpha_x g}_\lxq \nm{\agl{q} \nabla_q \nabla^\alpha_x g}_\lxq \no\\
	  \le\ & \tfrac{1}{10} \nm{\agl{q} \nabla_q \nabla^\alpha_x g}^2_\lxq + C \nm{\nabla^\alpha_x g}^2_\lxq,
	\end{align}
from which it follows that, for the summation over $|\alpha| \le 3$,
	\begin{align}
		\sum_{|\alpha|\le 3}\skpm{\L \nabla^\alpha_x g}{\agl{q}^2 \nabla^\alpha_x g}
		\ge \tfrac{9}{10}\nm{\agl{q} \nabla_q g}^2_\hxq{3} - C \nm{g}^2_\hxq{3}.
	\end{align}

In what follows, we turn to estimate the right-hand side of equation \eqref{eq:deriv g}. By straightforward calculations, we get
	\begin{align}
	  \skpm{\nabla^\alpha_x \div u}{\agl{q}^2 \nabla^\alpha_x g}
	  \le\ & |\nabla^\alpha_x \div u|_{L^2_x} \left( \int \agl{q}^2 M \d q \right)^{1/2} \nm{\agl{q} \nabla^\alpha_x g}_\lxq \\\no
	  \ls\ & |\nabla_x u|_\hx{3} \nm{\agl{q} \nabla^\alpha_x g}_\lxq.
	\end{align}
Next we consider the following term
	\begin{align}
	 R \triangleq \skpm{\nabla^\alpha_x (\div u \cdot g)}{\agl{q}^2 \nabla^\alpha_x g}
	 =\ & \sum_{|\alpha|=|\alpha_1|+|\alpha_2|} \skpm{\nabla^{\alpha_1}_x \div u \nabla^{\alpha_2}_x g}{\agl{q}^2 \nabla^\alpha_x g} \\\no
	  \triangleq\ & \sum_{|\alpha|=|\alpha_1|+|\alpha_2|} R_{|\alpha_1|,|\alpha_2|},
	\end{align}
which can be estimated by a similar process of discussing the values of $\alpha_1$ and $\alpha_2$, more precisely, we have
	\begin{align}
	  R_{3,0} \le\ & |\nabla^3_x \div u|_\lx \nm{\agl{q} g}_{L^\infty_x L^2_{M\!q}} \nm{\agl{q} \nabla^\alpha_x g}_\lxq \\\no
	  \ls\ & |\nabla_x u|_\hx{3} \nm{\agl{q} g}_\hxq{2} \nm{\agl{q} \nabla^\alpha_x g}_\lxq, \\\no
	  R_{2,1}
	  \ls\ & |u|_\hx{3} \nm{\agl{q} \nabla^1_x g}_\hxq{2} \nm{\agl{q} \nabla^\alpha_x g}_\lxq, \\\no
	  R_{1,2}
	  \ls\ & |u|_\hx{3} \nm{\agl{q} \nabla^2_x g}_\hxq{1} \nm{\agl{q} \nabla^\alpha_x g}_\lxq, \\\no
	  R_{0,3}
	  \ls\ & |u|_\hx{3} \nm{\agl{q} \nabla^\alpha_x g}^2_\lxq,
	\end{align}
combining these above inequalities together gives that, for $|\alpha| \le 3$,
	\begin{align}
	  R = \skpm{\nabla^\alpha_x (\div u \cdot g)}{\agl{q}^2 \nabla^\alpha_x g}
	  \ls (|u|_\hx{3} + |\nabla_x u|_\hx{3}) \nm{\agl{q} g}^2_\hxq{3}.
	\end{align}
	
Now we deal with the last term involving the factor $(1+g) \nabla_x u q \nabla_q U$. Firstly, we notice there exists cancellation formula between equations \eqref{eq:deriv u} and \eqref{eq:deriv g}, i.e.
	\begin{align}\label{eq:cancellation}
	  \skp{\div \nabla^\alpha_x \int \nabla_q U \otimes q g M \d q}{\nabla^\alpha_x u}
	  + \skpm{\nabla^\alpha_x (\nabla_x u q \nabla_q U)}{\nabla^\alpha_x g} =0.
	\end{align}
Moreover, it is an easy matter to get
	\begin{align}
		\skpm{\nabla^\alpha_x (\nabla_x u q \nabla_q U)}{\agl{q}^2 \nabla^\alpha_x g}
		=\ & \iint (\nabla^\alpha_x \nabla_x u) q \nabla_q U \agl{q}^2 \nabla^\alpha_x g M\d q\d x \\\no
		\le\ & |\nabla_x u|_\hx{3} \nm{q \nabla_q U \nabla^\alpha_x g}_\lxq (\int \agl{q}^4 M\d q)^{1/2} \\\no
		\ls\ & |\nabla_x u|_\hx{3} \nm{\agl{q} \nabla_q \nabla^\alpha_x g}_\lxq,
	\end{align}
where we have used Lemma \ref{lemm:weight inequ} again.

Now we are left to deal with the term
	\begin{align}
	 	\widetilde{R} \triangleq \skpm{\nabla^\alpha_x (\nabla_x u q \nabla_q U g)}{\agl{q}^2 \nabla^\alpha_x g}
	 =\ & \sum_{|\alpha|=|\alpha_1|+|\alpha_2|} \skpm{\nabla^{\alpha_1}_x \nabla_x u q \nabla_q U \nabla^{\alpha_2}_x g}{\agl{q}^2 \nabla^\alpha_x g} \no\\
	  \triangleq\ & \sum_{|\alpha|=|\alpha_1|+|\alpha_2|} \widetilde{R}_{|\alpha_1|,|\alpha_2|}.
	\end{align}
A discussing process over the values of $\alpha_1$ and $\alpha_2$ as above gives that
	\begin{align}
	  \widetilde{R}_{3,0}
	  \ls\ & |\nabla_x u|_\hx{3} \nm{\agl{q} \nabla_q g}_\hxq{2} \nm{\agl{q} \nabla_q \nabla^\alpha_x g}_\lxq, \\\no
	  \widetilde{R}_{2,1}
	  \ls\ & |u|_\hx{3} \nm{\agl{q} \nabla_q \nabla^1_x g}_\hxq{2} \nm{\agl{q} \nabla_q \nabla^\alpha_x g}_\lxq, \\\no
	  \widetilde{R}_{1,2}
	  \ls\ & |u|_\hx{3} \nm{\agl{q} \nabla_q \nabla^2_x g}_\hxq{1} \nm{\agl{q} \nabla_q \nabla^\alpha_x g}_\lxq, \\\no
	  \widetilde{R}_{0,3}
	  \ls\ & |u|_\hx{3} \nm{\agl{q} \nabla_q \nabla^\alpha_x g}^2_\lxq,
	\end{align}
then from which we can infer that, for $|\alpha| \le 3$,
	\begin{align}
	  \widetilde{R} =\ & \skpm{\nabla^\alpha_x (\nabla_x u q \nabla_q U g)}{\agl{q}^2 \nabla^\alpha_x g} \\\no
	  	\ls\ & |\nabla_x u|_\hx{3} \nm{\agl{q} \nabla_q g}_\hxq{2} \nm{\agl{q} \nabla_q g}_\hxq{3} + |u|_\hx{3} \nm{\agl{q} \nabla_q g}^2_\hxq{3}.
	\end{align}

Finally, we add up all the above inequalities and take summation over $|\alpha| \le 3$ to get the higher order spatial derivatives estimates that
	\begin{align}\label{esm:deriv g-x}
	  & \tfrac{1}{2} \tfrac{\d}{\d t} \nm{\agl{q} g}^2_\hxq{3} + \tfrac{9}{10} \nm{\agl{q} \nabla_q g}^2_\hxq{3} \\\no
	  \ls\ & (|u|_\hx{3} + |\nabla_x u|_\hx{3}) \nm{\agl{q} g}^2_\hxq{3} + \nm{g}^2_\hxq{3}
	  			+ |u|_\hx{3} \nm{\agl{q} \nabla_q g}^2_\hxq{3} \\\no
	  		& + |\nabla_x u|_\hx{3} \nm{\agl{q} \nabla_q g}_\hxq{2} \nm{\agl{q} \nabla_q g}_\hxq{3}
	  			+ |\nabla_x u|_\hx{3} \nm{\agl{q} \nabla_q g}_\hxq{3} .
	\end{align}
As a straightforward corollary result, it is easy to check
	\begin{align}\label{esm:deriv g-x pure}
	  & \tfrac{1}{2} \tfrac{\d}{\d t} \nm{g}^2_\hxq{3} + \nm{\nabla_q g}^2_\hxq{3} \\\no
	  \ls\ & (|u|_\hx{3} + |\nabla_x u|_\hx{3}) \nm{g}^2_\hxq{3}
	  			+ |u|_\hx{3} \nm{\agl{q} \nabla_q g}_\hxq{3} \nm{g}_\hxq{3} \\\no
	  		& + |\nabla_x u|_\hx{3} \nm{\agl{q} \nabla_q g}_\hxq{2} \nm{g}_\hxq{3}
	  		+ \sum_{|\alpha|\le 3} \skpm{\nabla^\alpha_x (\nabla_x u q \nabla_q U)}{\nabla^\alpha_x g} .
	\end{align}

%% ---------------------------------------------------------------------
\subsubsection{Estimates for mixed derivatives} % (fold)
\label{ssub:estimates_for_mixed_derivatives}

% subsubsection estimates_for_mixed_derivatives (end)

%% ---------------------------------------------------------------------

Applying $\nabla^\alpha_\beta$ with the restrictions $|\alpha| + |\beta|\le 3 $ and $|\beta| \ge 1$ to equation \eqref{eq:C-polymeric-pert}$_3$, we get
	\begin{align}\label{eq:deriv g mixed}
	  \p_t \nabla^\alpha_\beta g + \nabla^\alpha_\beta (u\cdot \nabla_x g) & + \nabla^\alpha_\beta (\nabla_x u q \nabla_q g) + \nabla^\alpha_\beta \mathcal{L} g \\\no
		 	& = - \nabla^\alpha_\beta [2(1+g) \div u + (1+g) \nabla_x u q \nabla_q U].
	\end{align}

Taking $\lxq$ inner product with the quantity $\agl{q}^2 \nabla^\alpha_\beta g$, we get firstly,
	\begin{align*}
	  \skpm{\p_t \nabla^\alpha_\beta g}{\agl{q}^2 \nabla^\alpha_\beta g}
	  = \tfrac{1}{2} \tfrac{\d}{\d t} \nm{\agl{q} \nabla^\alpha_\beta g}^2_\lxq.
	\end{align*}

In the following texts, we will use frequently the method of discussing by cases. For the sake of exposition, we introduce
	\begin{align*}
		\Lambda \triangleq\ & \left\{ |\alpha|=|\alpha_1|+|\alpha_2|,\ |\beta|\ge 1 \right\}, \\
	  \Lambda' \triangleq\ & \left\{ |\alpha|=|\alpha_1|+|\alpha_2|,\ |\beta|=|\beta_1|+|\beta_2|,\ |\beta|\ge 1 \right\},
	\end{align*}
where $|\alpha| + |\beta| = 3$. It should be pointed out that, we will split the set into three cases, corresponding to three values of $|\alpha_1|=0,\ 1,\ 2$, separately. Then it is easy to infer that, the triple $(|\alpha_1|,\ |\alpha_2|,\ |\beta|)$ will take values of $(0,\ |\alpha|,\ |\beta|)$, $(1,\ |\alpha|-1,\ |\beta|)$ and $(2,\ 0,\ 1)$, separately.

As for the convection term, we denote that
	\begin{align}
	 W \triangleq \skpm{\nabla^\alpha_\beta (u\cdot \nabla_x g)}{\agl{q}^2 \nabla^\alpha_\beta g}
	 =\ & \sum_{\Lambda} \skpm{\nabla^{\alpha_1}_x u \nabla_x \nabla^{\alpha_2}_\beta g}{\agl{q}^2 \nabla^\alpha_ \beta g} \no\\
	  \triangleq\ & \sum_{\Lambda} W_{|\alpha_1|,|\alpha_2|}.
	\end{align}
We then get separately that
	\begin{align}
		W_{0,|\alpha|} \le\ & |\div u|_{L^\infty_x} \nm{\agl{q} \nabla^\alpha_\beta g}^2_\lxq
		\ls |u|_\hx{3} \nm{\agl{q} \nabla^\alpha_\beta g}^2_\lxq, \\\no
		W_{1,|\alpha|-1} \le\ & |\nabla^1_x u|_{L^\infty_x} \nm{\agl{q} \nabla^\alpha_\beta g}^2_\lxq
		\ls |u|_\hx{3} \nm{\agl{q} \nabla^\alpha_\beta g}^2_\lxq, \\\no
		W_{2,0} \le\ & |\nabla^2_x u|_{L^4_x} \nm{\agl{q} \nabla^1_1 g}_{L^4_x L^2_{M\!q}} \nm{\agl{q} \nabla^\alpha_\beta g}_\lxq \\\no
		\ls\ & |u|_\hx{3} \nm{\agl{q} \nabla^1_1 g}_\hxq{1} \nm{\agl{q} \nabla^\alpha_\beta g}_\lxq,
	\end{align}
which together gives that, for $|\alpha| + |\beta| \le 3$,
	\begin{align}
	  W \ls |u|_\hx{3} \nm{\agl{q} \nabla^\beta_q g}_\hxq{{|\alpha|}} \nm{\agl{q} \nabla^\alpha_\beta g}_\lxq.
	\end{align}
	
For the third term on the left-hand side of equation \eqref{eq:deriv g mixed}, we have
	\begin{align}
	  \Gamma \triangleq\ & \skpm{\nabla^\alpha_\beta (\nabla_x u q \nabla_q g)}{\agl{q}^2 \nabla^\alpha_\beta g} \\\no
	 =\ & \sum_{\Lambda} \skpm{\nabla^{\alpha_1 +1}_x u q \nabla_q \nabla^{\alpha_2}_\beta g}{\agl{q}^2 \nabla^\alpha_ \beta g}
	 			+ \sum_{\Lambda} \skpm{\nabla^{\alpha_1 +1}_x u \nabla_q \nabla^{\alpha_2}_{\beta-1} g \agl{q}}{\agl{q} \nabla^\alpha_ \beta g} \\\no
	  \triangleq\ & \sum_{\Lambda} (\Gamma^{1}_{|\alpha_1|,|\alpha_2|} + \Gamma^{2}_{|\alpha_1|,|\alpha_2|}).
	\end{align}

Estimate these above terms by using the Young inequality, the Sobolev embedding inequality and Lemma \ref{lemm:weight inequ} as we did before, then we list them as follows,
	\begin{align}
		\Gamma^{1}_{0,|\alpha|} \le\ & |\nabla_x u|_{L^\infty_x} \nm{|q| \nabla_q \nabla^\alpha_\beta g}_\lxq \nm{\agl{q}^2 \nabla^\alpha_\beta g}_\lxq \\\no
			\ls\ & |u|_\hx{3} \nm{\agl{q} \nabla_q \nabla^\alpha_\beta g}_\lxq (\nm{\agl{q} \nabla_q \nabla^\alpha_\beta g}_\lxq + \nm{\nabla^\alpha_\beta g}_\lxq) , \\\no
		\Gamma^{1}_{1,|\alpha|-1} \le\ & |\nabla^2_x u|_{L^4_x} \nm{|q| \nabla_q \nabla^{\alpha-1}_\beta g}_{L^4_x L^2_{M\!q}} \nm{\agl{q}^2 \nabla^\alpha_\beta g}_\lxq \\\no
			\ls\ & |u|_\hx{3} \nm{\agl{q} \nabla_q \nabla^{\alpha-1}_\beta g}_\hxq{1} (\nm{\agl{q} \nabla_q \nabla^\alpha_\beta g}_\lxq + \nm{\nabla^\alpha_\beta g}_\lxq) , \\\no
		\Gamma^{1}_{2,0} \le\ & |\nabla^3_x u|_{L^2_x} \nm{|q| \nabla_q \nabla^0_1 g}_{L^\infty_x L^2_{M\!q}} \nm{\agl{q}^2 \nabla^\alpha_\beta g}_\lxq \\\no
			\ls\ & |u|_\hx{3} \nm{\agl{q} \nabla_q \nabla^0_1 g}_\hxq{2} (\nm{\agl{q} \nabla_q \nabla^\alpha_\beta g}_\lxq + \nm{\nabla^\alpha_\beta g}_\lxq) ,
	\end{align}
and similarly,
	\begin{align}
		\Gamma^{2}_{0,|\alpha|}
		\ls\ & |u|_\hx{3} \nm{\agl{q} \nabla^\alpha_\beta g}^2_\lxq, \\\no
		\Gamma^{2}_{1,|\alpha|-1}
		\ls\ & |u|_\hx{3} \nm{\agl{q} \nabla^{\alpha-1}_\beta g}_\hxq{1} \nm{\agl{q} \nabla^\alpha_\beta g}_\lxq, \\\no
		\Gamma^{2}_{2,0}
		\ls\ & |u|_\hx{3} \nm{\agl{q} \nabla^0_1 g}_\hxq{2} \nm{\agl{q} \nabla^\alpha_\beta g}_\lxq.
	\end{align}
Therefore we can deduce that, for $|\alpha| + |\beta| \le 3$,
	\begin{align}
	  \Gamma \ls\ & |u|_\hx{3} \nm{\agl{q} \nabla_q \nabla^\beta_q g}_\hxq{{|\alpha|}} \left( \nm{\agl{q} \nabla_q \nabla^\alpha_\beta g}_\lxq + \nm{\nabla^\alpha_\beta g}_\lxq \right) \\\no
	  & + |u|_\hx{3} \nm{\agl{q} \nabla^\beta_q g}_\hxq{{|\alpha|}} \nm{\agl{q} \nabla^\alpha_\beta g}_\lxq.
	\end{align}

Notice the following elementary fact
	\begin{align}
	  \nabla^\alpha_ \beta \L g = - \nabla^\alpha_ \beta (\Delta_q g - \nabla_q U \nabla_q g)
	  =\ & - (\Delta_q \nabla^\alpha_ \beta g - \nabla_q U \nabla_q \nabla^\alpha_ \beta g)
	  	+ [\nabla^\alpha_ \beta, \nabla_q U] \nabla_q g \\\no
	 	=\ & \L \nabla^\alpha_ \beta g + [\nabla^\alpha_ \beta, \nabla_q U] \nabla_q g,
	\end{align}
then the micro-dissipation contributions can be rewritten as
	\begin{align}\label{temp}
	  \skpm{\nabla^\alpha_ \beta \L g}{\agl{q}^2 \nabla^\alpha_ \beta g}
	  =\ & \skpm{\L \nabla^\alpha_ \beta  g}{\agl{q}^2 \nabla^\alpha_ \beta g}
	  		+ \skpm{[\nabla^\alpha_ \beta, \nabla_q U] \nabla_q g}{\agl{q}^2 \nabla^\alpha_ \beta g} \\\no
	  =\ & \skpm{\nabla_q \nabla^\alpha_ \beta g}{\agl{q}^2 \nabla_q \nabla^\alpha_ \beta g}
	  			+ \skpm{\nabla_q \nabla^\alpha_ \beta g}{2 q \nabla^\alpha_ \beta g} \\\no
	   	 & + \skpm{[\nabla^\alpha_ \beta, \nabla_q U] \nabla_q g}{\agl{q}^2 \nabla^\alpha_\beta g}.
 	\end{align}
Together with Lemma \ref{lemm:weight inequ}, the assumptions on the potential $U$ \eqref{asmp-2} yields the commutator estimates
 \begin{align}\label{temp1}
   \skpm{[\nabla^\alpha_ \beta, \nabla_q U] \nabla_q g}{\agl{q}^2 \nabla^\alpha_\beta g}
   =\ & \sum_{\substack{|\beta|= |\beta_1|+|\beta_2|\\|\beta_1| \ge 1,\ |\beta_2| \le |\beta|-1}}
   		\skpm{\agl{q} \nabla^{\beta_1}_q (\nabla_q U) \nabla^\alpha_{\beta_2+1} g}{\agl{q} \nabla^\alpha_\beta g} \\\no
   \ls\ & \sum_{1 \le |\beta'| \le |\beta|} \iint \agl{q \nabla_q U} \nabla^\alpha_\beta g \cdot \agl{q} \nabla^\alpha_{\beta'} g M\d q\d x \\\no
   \ls\ & \sum_{1 \le |\beta'| \le |\beta|} \nm{\agl{q} \nabla^\alpha_{\beta'} g}_\lxq (\nm{\agl{q} \nabla_q \nabla^\alpha_\beta g}_\lxq + \nm{\nabla^\alpha_\beta g}_\lxq).
 \end{align}

On the other hand, combining with the H\"{o}lder inequality and Lemma \ref{lemm:weight inequ} together yields that
	\begin{align}\label{temp2}
		\skpm{\nabla_q \nabla^\alpha_ \beta g}{2 q \nabla^\alpha_ \beta g}
		=\ & 3 \nm{\nabla^\alpha_ \beta g}^2_\lxq - \iint q \nabla_q U |\nabla^\alpha_ \beta g|^2 M\d q\d x \\\no
	  \ls\ & 3 \nm{\nabla^\alpha_ \beta g}^2_\lxq + \nm{\nabla^\alpha_ \beta g}_\lxq \nm{q \nabla_q U \nabla^\alpha_ \beta g}_\lxq \\\no
	  \ls\ & C \nm{\nabla^\alpha_ \beta g}^2_\lxq + \nm{\nabla^\alpha_ \beta	g}_\lxq \nm{\agl{q} \nabla_q \nabla^\alpha_ \beta g}_\lxq \no\\
	  \le\ & \tfrac{1}{10} \nm{\agl{q} \nabla_q \nabla^\alpha_ \beta g}^2_\lxq + C \nm{\nabla^\alpha_ \beta g}^2_\lxq.
	\end{align}

Thus, we can infer from inserting equations \eqref{temp1}-\eqref{temp2} into \eqref{temp} that,
	\begin{align}
		& \skpm{\nabla^\alpha_ \beta \L g}{\agl{q}^2 \nabla^\alpha_ \beta g} \\\no
		\ge\ & \tfrac{9}{10}\nm{\agl{q} \nabla^\alpha_ \beta \nabla_q g}^2_\lxq
			  - C \sum_{1 \le |\beta'| \le |\beta|} \nm{\agl{q} \nabla^\alpha_{\beta'} g}_\lxq \left( \nm{\agl{q} \nabla_q \nabla^\alpha_\beta g}_\lxq + \nm{\nabla^\alpha_\beta g}_\lxq \right) .
	\end{align}

We turn to deal with the right-hand side terms of equation \eqref{eq:deriv g mixed}. Firstly, notice that the fact
	\begin{align}
	  \skpm{\nabla^\alpha_ \beta \div u}{\agl{q}^2 \nabla^\alpha_ \beta g} = 0,
	\end{align}
then it requires to consider
	\begin{align}
	  \Theta = \skpm{\nabla^\alpha_ \beta (\div u g)}{\agl{q}^2 \nabla^\alpha_ \beta g}
	  =\ & \sum_{\Lambda} \skpm{\nabla^{\alpha_1} \div u \nabla^{\alpha_2}_\beta g}{\agl{q}^2 \nabla^\alpha_ \beta g} \\\no
	  =\ & \sum_{\Lambda} \Theta_{|\alpha_1|,|\alpha_2|},
	\end{align}
where the three terms $\Theta_{|\alpha_1|,|\alpha_2|}$ can be estimated similarly as before. For simplicity, we only write the result:
	\begin{align}
	  \Theta
	  \ls |u|_\hx{3} \nm{\agl{q} \nabla^\beta_q g}_\hxq{{|\alpha|}} \nm{\agl{q} \nabla^\alpha_\beta g}_\lxq.
	\end{align}
	
By the assumptions on $U$ \eqref{asmp-2} and Lemma \ref{lemm:weight inequ}, it is easy to get
	\begin{align}
		\skpm{\nabla^\alpha_ \beta (\nabla_x u q \nabla_q U)}{\agl{q}^2 \nabla^\alpha_ \beta g}
		=\ & \iint (\nabla^\alpha_x \nabla_x u) \nabla^\beta_q (q \nabla_q U) \agl{q}^2 \nabla^\alpha_ \beta g M\d q\d x \\\no
		\le\ & |\nabla^\alpha_x \nabla_x u|_\lx \nm{\agl{q \nabla_q U} \nabla^\alpha_ \beta g}_\lxq (\int \agl{q}^4 M\d q)^{1/2} \\\no
		\ls\ & |u|_\hx{3} (\nm{\agl{q} \nabla_q \nabla^\alpha_ \beta g}_\lxq + \nm{\nabla^\alpha_ \beta g}_\lxq).
	\end{align}

At last, we need to deal with the term
	\begin{align}
	 	\Xi \triangleq\ & \skpm{\nabla^\alpha_ \beta (\nabla_x u q \nabla_q U g)}{\agl{q}^2 \nabla^\alpha_ \beta g} \\\no
	 =\ & \sum_{\Lambda'} \skpm{\nabla^{\alpha_1}_x \nabla_x u \nabla^{\beta_1}_q (q \nabla_q U) \nabla^{\alpha_2}_{\beta_2} g}{\agl{q}^2 \nabla^\alpha_\beta g}
	  \triangleq \sum_{\Lambda'} \Xi_{|\alpha_1|,|\alpha_2|}.
	\end{align}
A discussing process over the values of $\alpha_1$ and $\alpha_2$ as above gives that
	\begin{align}
	  \Xi_{0,|\alpha|}
	  	\le & \sum_{|\beta_2|\le |\beta|} |\nabla_x u|_{L^\infty_x} \nm{\agl{q \nabla_q U} \nabla^\alpha_{\beta_2} g}_\lxq \nm{\agl{q}^2 \nabla^\alpha_ \beta g}_\lxq \\\no
	  	\ls & \sum_{|\beta_2|\le |\beta|} |u|_\hx{3} (\nm{\agl{q} \nabla_q \nabla^\alpha_{\beta_2} g}_\lxq + \nm{\nabla^\alpha_{\beta_2} g}_\lxq) (\nm{\agl{q} \nabla_q \nabla^\alpha_ \beta g}_\lxq + \nm{\nabla^\alpha_ \beta g}_\lxq), \\\no
	  \Xi_{1,|\alpha|-1}
	  	\le & \sum_{|\beta_2|\le |\beta|} |\nabla^1_x \nabla_x u|_{L^4_x} \nm{\agl{q \nabla_q U} \nabla^{\alpha-1}_{\beta_2} g}_{L^4_x L^2_{M\!q}} \nm{\agl{q}^2 \nabla^\alpha_ \beta g}_\lxq \\\no
	  	\ls & \sum_{|\beta_2|\le |\beta|}\!\!\!\! |u|_\hx{3} (\nm{\agl{q} \!\! \nabla_q \nabla^{\alpha-1}_{\beta_2} g}_\hxq{1} + \nm{\nabla^{\alpha-1}_{\beta_2} g}_\hxq{1}) (\nm{\agl{q} \!\! \nabla_q \nabla^\alpha_ \beta g}_\lxq + \nm{\nabla^\alpha_ \beta g}_\lxq), \\\no
	  \Xi_{2,0}
	  	\le & \sum_{|\beta_2|\le |\beta|} |\nabla^1_x \nabla_x u|_{L^2_x} \nm{\agl{q \nabla_q U} \nabla^0_{\beta_2} g}_{L^\infty_x L^2_{M\!q}} \nm{\agl{q}^2 \nabla^\alpha_ \beta g}_\lxq \\\no
	  	\ls & \sum_{|\beta_2|\le |\beta|} |u|_\hx{3} (\nm{\agl{q} \nabla_q \nabla^0_{\beta_2} g}_\hxq{2} + \nm{\nabla^0_{\beta_2} g}_\hxq{2}) (\nm{\agl{q} \nabla_q \nabla^\alpha_ \beta g}_\lxq + \nm{\nabla^\alpha_ \beta g}_\lxq),
	\end{align}
then from which we can infer that,
	\begin{align}
	  \Xi \ls \sum_{|\beta_2|\le |\beta|} |u|_\hx{3}
	  		(\nm{\agl{q} \nabla_q \nabla^{\beta_2}_q g}_\hxq{{|\alpha|}} + \nm{\nabla^{\beta_2}_q g}_\hxq{{|\alpha|}})
	  		(\nm{\agl{q} \nabla_q \nabla^\alpha_\beta g}_\lxq + \nm{\nabla^\alpha_\beta g}_\lxq) .
	\end{align}

At the end, combining all the above inequalities leads us to the mixed derivatives estimates that
	\begin{align}\label{esm:deriv g-q}
	  & \tfrac{1}{2} \tfrac{\d}{\d t} \nm{\agl{q} \nabla^\alpha_ \beta g}^2_\lxq + \tfrac{9}{10} \nm{\agl{q} \nabla_q \nabla^\alpha_ \beta g}^2_\lxq \\\no
	 \ls\ & |u|_\hx{3} \nm{\agl{q} \nabla^\beta_q g}_\hxq{{|\alpha|}} \nm{\agl{q} \nabla^\alpha_\beta g}_\lxq
	  		+ |u|_\hx{3} (\nm{\agl{q} \nabla_q \nabla^\alpha_ \beta g}_\lxq + \nm{\nabla^\alpha_ \beta g}_\lxq) \\\no
	  		& + |u|_\hx{3} \nm{\agl{q} \nabla_q \nabla^\beta_q g}_\hxq{{|\alpha|}}
	  			( \nm{\agl{q} \nabla_q \nabla^\alpha_\beta g}_\lxq + \nm{\nabla^\alpha_\beta g}_\lxq ) \\\no
	  		& + \sum_{1 \le |\beta'| \le |\beta|} \nm{\agl{q} \nabla^\alpha_{\beta'} g}_\lxq (\nm{\agl{q} \nabla_q \nabla^\alpha_\beta g}_\lxq + \nm{\nabla^\alpha_\beta g}_\lxq) \\\no
	  		& + \sum_{|\beta_2|\le |\beta|} |u|_\hx{3}
	  		(\nm{\agl{q} \nabla_q \nabla^{\beta_2}_q g}_\hxq{{|\alpha|}} + \nm{\nabla^{\beta_2}_q g}_\hxq{{|\alpha|}})
	  		(\nm{\agl{q} \nabla_q \nabla^\alpha_\beta g}_\lxq + \nm{\nabla^\alpha_\beta g}_\lxq).
	\end{align}

%% ---------------------------------------------------------------------
%% ---------------------------------------------------------------------

\subsection{The a priori estimates for the compressible polymeric system} % (fold)
\label{sub:the_a_priori_estimates_for_the_compressible_polymeric_systems}

% subsection the_a_priori_estimates_for_the_compressible_polymeric_systems (end)

%% ---------------------------------------------------------------------

Now we are in the position to combine all the above higher order derivatives estimates to get the a priori estimates for the compressible polymeric system \eqref{eq:C-polymeric-pert}. We first introduce a variation of the mixed derivative norm $\nm{f}^2_\hxqm{s}$. More precisely, for some sufficiently small (and fixed) constant $\eta >0$, we define
	\begin{align}\label{def:small eta}
	  \nm{\agl{q} f}^2_\hxqe{s} = \sum_{|\alpha| \le s} \eta \nm{\agl{q} \nabla^\alpha_x f|}^2_\lxq + \sum_{|\alpha|+|\beta|\le s,\ |\beta| \ge 1} \eta^{|\beta|} \nm{\agl{q} \nabla^\alpha_\beta f|}^2_\lxq.
	\end{align}
Note that $\nm{\cdot}_\hxqe{0} \sim \nm{\cdot}_\hxqm{0} = \nm{\cdot}_\lxq$ and $\nm{\cdot}_\hxqe{s} \sim \nm{\cdot}_\hxqm{s}$.

By this definition, it follows from combining equations \eqref{esm:deriv g-x} and \eqref{esm:deriv g-q} for all derivatives $|\alpha|+|\beta| \le 3$, that
	\begin{align}\label{esm:micro-eta}
	  & \tfrac{1}{2} \tfrac{\d}{\d t} \nm{\agl{q} g}^2_\hxqe{3} + \tfrac{9}{10} \nm{\agl{q} \nabla_q g}^2_\hxqe{3} \\\no
	  \ls\ & |u|_\hx{3} \nm{\agl{q} \nabla_q g}^2_\hxqe{3}
	  		+ \eta^{\frac{1}{2}} |u|_\hx{3} \nm{\agl{q} \nabla_q g}_\hxqe{3}
	  		+ \eta^{\frac{1}{2}} \nm{\agl{q} \nabla_q g}^2_\hxqe{3}
	  		+ \eta \nm{\nabla_q g}^2_\hxq{3} \\\no
 	  		& + |\nabla_x u|_\hx{3} \nm{\agl{q} \nabla_q g}_\hxqe{3} \nm{\agl{q} g}_\hxqe{3}
 	  		+ \eta^{\frac{1}{2}} |\nabla_x u|_\hx{3} \nm{\agl{q} \nabla_q g}_\hxqe{3} . 			
	\end{align}

Denote that
	\begin{align}\label{eq:energy-eta}
	  E_\eta =\ & |\rho|^2_\hx{3} + |u|^2_\hx{3} + \nm{g}^2_\hxq{3} + \nm{\agl{q} g}^2_\hxqe{3} + \eta \Sigma_{|\alpha| \le 2} \skp{\nabla^\alpha_x u}{\nabla_x \nabla^\alpha_x \rho}, \\\no
	  D_\eta =\ & \mu|\nabla_x u|^2_\hx{3} + (\mu+\xi)|\div u|^2_\hx{3} + \nm{\nabla_q g}^2_\hxq{3} + \nm{\agl{q} \nabla_q g}^2_\hxqe{3} + \eta |\nabla_x \rho|^2_\hx{2}.
	\end{align}
It is obvious that $\eta^{3} E(t) \le E_\eta(t) \le 2 E(t)$ and $\eta^{3} D(t) \le D_\eta(t)$, hence the small assumption $E(t) \le \ep$ yields that $E_\eta(t) \ls \ep$.

Then by combining the above estimates \eqref{esm:fluid}, \eqref{esm:density dispa}, \eqref{esm:deriv g-x pure} and \eqref{esm:micro-eta}, we can infer finally that,
	\begin{align}\label{esm:total}
	  \tfrac{1}{2} \tfrac{\d}{\d t} E_\eta (t)
	  	+ \tfrac{4}{5} D_\eta (t)
	  \ls (\ep^{\frac{1}{2}} + \eta^{\frac{1}{2}}) D_\eta(t).
	\end{align}	
As a consequence, choosing some sufficiently small $\eta>0$ immediately yields that, for any $0 <t < T$,
	\begin{align}
	  E_\eta(t) + \int_0^t D_\eta (t) \d t \le E_\eta(0),
	\end{align}
which concludes the a priori estimates stated in Proposition \ref{prop:A-priori} using again the equivalence of $E_\eta(0) \sim E(0)$, more precisely, there exists some constant $C_0>1$ such that,
	\begin{align}
	  E (t) \le C_0 E(0).
	\end{align}

%% ------------------------------------------------------------------------ %%

\section{Global existence with small initial assumptions} % (fold)
\label{sec:global_existence_with_small_initial_assumptions}

% section global_existence_with_small_initial_assumptions (end)

In this section, we aim mainly at justifying the global existence of classic solutions to the compressible polymeric system \eqref{eq:C-polymeric-pert}, under some small fluctuational assumptions near the global equilibrium state $(1,\, 0,\, M)$.

The proof of existence proceeds as follows: First, by constructing a sequence of approximate solution to the system \eqref{eq:C-polymeric-pert} using a standard iteration scheme, we can prove the uniform-in-$n$ estimate of the approximate sequence $(\rho^n,u^n,g^n)$; second, we prove that the approximate sequence is convergent in some low norm. Then performing one standard scheme involving interpolation theory and Fatou lemma will gives the local-in-time existence and uniqueness result of solution to the original system \eqref{eq:C-polymeric-pert}, based on which we can finally obtain the global-in-time existence by a continuum argument.

More precisely, we construct the iterating approximating sequence as follows,
	\begin{align}\label{eq:appro}
	  \begin{cases}
	  	\p_t \rho^{n+1} + u^n \nabla_x \rho^{n+1} + (1+\rho^n) \div u^{n+1} = 0,
		\\[5pt]
			\begin{aligned}
				\p_t u^{n+1} + u^n \cdot \nabla_x u^{n+1} + & \tfrac{P'(1+\rho^n)}{1+\rho^n} \nabla_x \rho^{n+1}
			\\[3pt]
				& = \tfrac{1}{1+\rho^n} \div \Sigma(u^{n+1}) + \tfrac{1}{1+\rho^n} \div \int_{\mathbb{R}^3} \nabla_q U \otimes q g^{n+1} M \d q ,
			\end{aligned}
	  \\[7pt]
			\begin{aligned}
				\p_t g^{n+1} + u^n \cdot \nabla_x g^{n+1} + & \nabla_x u^n q \nabla_q g^{n+1} + \mathcal{L} g^{n+1}
			\\[3pt]
				& = - 2 (1+g^{n+1}) \div u^n + (1+g^{n+1}) \nabla_x u^n q \nabla_q U ,
			\end{aligned}	
	  \end{cases}
	\end{align}
endowed with the initial data $(\rho^{n+1},u^{n+1},g^{n+1})|_{t=0}=(\rho_0(x), u_0(x), g_0(x,q))$. And moreover, we start with $(\rho^{0}(t,x),u^{0}(t,x),g^{0}(t,x,q)) \equiv (\rho_0(x), u_0(x), g_0(x,q))$.

\subsection{Uniform bound in a large norm} % (fold)
\label{sub:uniform_bound_in_a_large_norm}

% subsection uniform_bound_in_a_large_norm (end)

We state here the uniform-in-$n$ boundedness of the approximate sequence $(\rho^n,u^n,g^n)$ in the following lemma, which plays an important role in the proof of global existence.
\begin{lemma}\label{lemm:uniform esm}
	There exist $M_0 >0$ and $T_*>0$, such that if $\E(0) \le M_0/2$ and $\sup_{t\in [0,T_*]} \E_n(t) \le M_0$, then
		\begin{align}
		  \sup_{t\in [0,T_*]} \E_{n+1}(t) \le M_0.
		\end{align}
		
\end{lemma}

Corresponding to the previous energy estimates \eqref{esm:deriv rho} and \eqref{esm:deriv u}, it is an easy matter to get that
	\begin{align}\label{esm:rho-sequ}
	   \tfrac{1}{2} \tfrac{\d}{\d t} |\rho^{n+1}|^2_\hx{3}
	  \ls\ & |u^n|_\hx{3} |\rho^{n+1}|^2_\hx{3}
	  		+ |\rho^n|_\hx{3} |\rho^{n+1}|_\hx{3} |u^{n+1}|_\hx{3} \\\no
	  		& + (1+|\rho^n|_\hx{3}) |\rho^{n+1}|_\hx{3} |\nabla_x u^{n+1}|_\hx{3} ,
	\end{align}
and
	\begin{align}\label{esm:u-sequ}
	  & \tfrac{1}{2} \tfrac{\d}{\d t} |u^{n+1}|^2_\hx{3} + \tfrac{4}{5} \left( \mu|\nabla_x u^{n+1}|^2_\hx{3} + (\mu+\xi)|\div u^{n+1}|^2_\hx{3} \right) \\\no
	  \ls \ & |u^n|_\hx{3} |u^{n+1}|^2_\hx{3}
	  				+ |\rho^n|_\hx{3} |\nabla_x u^{n+1}|_\hx{3}
	  			  + |\rho^n|_\hx{3} |\rho^{n+1}|_\hx{3} |u^{n+1}|_\hx{3}
	  			\\\no
	  			& + |\rho^n|_\hx{3} |u^{n+1}|_\hx{3} |\nabla_x u^{n+1}|_\hx{3}
	  			  + (1+ |\rho^n|_\hx{3}) |u^{n+1}|_\hx{3} \nm{\nabla_q g^{n+1}}_\hxq{3}.
	\end{align}

On the other hand, performing almost the same process as that of equations\eqref{esm:deriv g-x} and \eqref{esm:deriv g-q},  we are able to infer from equation\eqref{eq:appro}$_3$ that,
	\begin{align}\label{esm:g-x sequ}
	  & \tfrac{1}{2} \tfrac{\d}{\d t} \nm{\agl{q} g^{n+1}}^2_\hxq{3} + \tfrac{9}{10} \nm{\agl{q} \nabla_q g^{n+1}}^2_\hxq{3} \\\no
	  \ls\ & (|u^n|_\hx{3} + |\nabla_x u^n|_\hx{3} +1) \nm{\agl{q} g^{n+1}}^2_\hxq{3}
	  			+ |u^n|_\hx{3} \nm{\agl{q} \nabla_q g^{n+1}}^2_\hxq{3} \\\no
	  		& + |\nabla_x u^n|_\hx{3} \nm{\agl{q} g^{n+1}}_\hxqm{3} \nm{\agl{q} \nabla_q g^{n+1}}_\hxq{3}
	  		  + |\nabla_x u^n|_\hx{3} \nm{\agl{q} g^{n+1}}_\hxq{3} .
	\end{align}
and
	\begin{align}\label{esm:g-q sequ}
	  & \tfrac{1}{2} \tfrac{\d}{\d t} \sum_{\Lambda} \nm{\agl{q} \nabla^\alpha_ \beta g^{n+1}}^2_\lxq + \tfrac{9}{10} \sum_{\Lambda} \nm{\agl{q} \nabla_q \nabla^\alpha_ \beta g^{n+1}}^2_\lxq \\\no
	  \ls\ & |u^n|_\hx{3} \left( \nm{\agl{q} g^{n+1}}^2_\hxqm{3} + \nm{\agl{q} \nabla_q g^{n+1}}^2_\hxqm{3} \right)
	  			+ \nm{\agl{q} g^{n+1}}^2_\hxqm{3} \\\no
	  		& + |u^n|_\hx{3} \nm{\agl{q} \nabla_q g^{n+1}}_\hxqm{3}
	  			+ \sum_{\Lambda} \sum_{1 \le |\beta'| \le |\beta|} \nm{\agl{q} \nabla^\alpha_{\beta'} g^{n+1}}_\lxq \nm{\agl{q} \nabla_q \nabla^\alpha_\beta g^{n+1}}_\lxq \\\no
	  		& + \sum_{\Lambda} \sum_{|\beta_2|\le |\beta|} |u^n|_\hx{3}
	  		\nm{\agl{q} \nabla_q \nabla^\alpha_{\beta_2} g^{n+1}}_\lxq \nm{\agl{q} \nabla_q \nabla^\alpha_\beta g^{n+1}}_\lxq \\\no
	  \ls\ & |u^n|_\hx{3} ( \nm{\agl{q} g^{n+1}}^2_\hxqm{3} + \nm{\agl{q} \nabla_q g^{n+1}}^2_\hxqm{3} )
	  			+ \nm{\agl{q} g^{n+1}}^2_\hxqm{3} \\\no
	  		& + (|u^n|_\hx{3} + \nm{\agl{q} g^{n+1}}_\hxqm{3}) \nm{\agl{q} \nabla_q g^{n+1}}_\hxqm{3}.
	\end{align}

Adding these above estimates entails that
	\begin{align}\label{esm:total-sequ}
	  & \tfrac{1}{2} \tfrac{\d}{\d t} (|\rho^{n+1}|^2_\hx{3} + |u^{n+1}|^2_\hx{3} + \nm{\agl{q} g^{n+1}}^2_\hxqm{3}) \\\no
	  	& + \tfrac{4}{5} \left( \mu|\nabla_x u^{n+1}|^2_\hx{3} + (\mu+\xi)|\div u^{n+1}|^2_\hx{3} + \nm{\agl{q} \nabla_q g^{n+1}}^2_\hxqm{3} \right) \\\no
	  \ls\ & |u^n|_\hx{3} (|\rho^{n+1}|^2_\hx{3} + |u^{n+1}|^2_\hx{3})
					 + |\rho^n|_\hx{3} |\rho^{n+1}|_\hx{3} |u^{n+1}|_\hx{3}
					 + (1+|\rho^n|_\hx{3}) |\rho^{n+1}|_\hx{3} |\nabla_x u^{n+1}|_\hx{3} \\\no
	  		 & + |\rho^n|_\hx{3} |u^{n+1}|_\hx{3} |\nabla_x u^{n+1}|_\hx{3}
	  		   + (1+|\rho^n|_\hx{3}) |u^{n+1}|_\hx{3} \nm{\nabla_q g^{n+1}}_\hxq{3} \\\no
	  		 & + (|u^n|_\hx{3} + |\nabla_x u^n|_\hx{3} + 1) \nm{\agl{q} g^{n+1}}^2_\hxqm{3}
	  		   + |u^n|_\hx{3} (\nm{\agl{q} \nabla_q g^{n+1}}_\hxqm{3} + \nm{\agl{q} \nabla_q g^{n+1}}^2_\hxqm{3}) \\\no
	  		 & + \nm{\agl{q} g^{n+1}}_\hxqm{3} (|\nabla_x u^n|_\hx{3} + \nm{\agl{q} \nabla_q g^{n+1}}_\hxqm{3}) \\\no
	  		 & + |\nabla_x u^n|_\hx{3} \nm{\agl{q} g^{n+1}}_\hxqm{3} \nm{\agl{q} \nabla_q g^{n+1}}_\hxqm{3} .
	\end{align}	
Let $E^n(t)= E(\rho^n,u^n,g^n)(t)$ and $D^n(t)= D(\rho^n,u^n,g^n)(t)$, then we can get the approximate energy estimate as follows:
	\begin{align}
	  & \tfrac{1}{2} \tfrac{\d}{\d t} E^{n+1}(t) + \tfrac{4}{5} D^{n+1}(t) \\\no
	  \le\ & C E^{n+1}(t) [(E^n)^{\frac{1}{2}}(t) + 1]
	  		   + C [1+(E^n)^{\frac{1}{2}}(t)] (E^{n+1})^{\frac{1}{2}}(t) (D^{n+1})^{\frac{1}{2}}(t)
	  		 	 + C E^{n+1}(t) (D^n)^{\frac{1}{2}}(t) \\\no
	  		 & + C (E^n)^{\frac{1}{2}}(t) [(D^{n+1})^{\frac{1}{2}}(t) + D^{n+1}(t)]
	  		 	 + C (E^{n+1})^{\frac{1}{2}}(t) (D^n)^{\frac{1}{2}}(t) \\\no
	  		 & + C (E^{n+1})^{\frac{1}{2}}(t) (D^n)^{\frac{1}{2}}(t) (D^{n+1})^{\frac{1}{2}}(t) .
	\end{align}
For the terms involving the factor $D^{n+1}$, we perform the H\"older inequalities to get, by noticing the induction assumptions $\sup_{t\in [0,T_*]} \E_n(t) \le M_0$, that
\begin{align*}
  (E^{n+1})^{\frac{1}{2}}(t) (D^{n+1})^{\frac{1}{2}}(t)
	  	\le\ & \tfrac{1}{15 C} D^{n+1}(t) + C E^{n+1}(t), \\
  (E^n)^{\frac{1}{2}}(t) (E^{n+1})^{\frac{1}{2}}(t) (D^{n+1})^{\frac{1}{2}}(t)
	  	\le\ & M_0^{\frac{1}{2}} D^{n+1}(t) + C M_0^{\frac{1}{2}} E^{n+1}(t), \\
  (E^n)^{\frac{1}{2}}(t) (D^{n+1})^{\frac{1}{2}}(t)
	  	\le\ & \tfrac{1}{15 C} D^{n+1}(t) + C M_0, \\
  (E^n)^{\frac{1}{2}}(t) D^{n+1}(t) \le\ & M_0^{\frac{1}{2}} D^{n+1}(t), \\
  (E^{n+1})^{\frac{1}{2}}(t) (D^n)^{\frac{1}{2}}(t) (D^{n+1})^{\frac{1}{2}}(t)
	  	\le\ & \tfrac{1}{15 C} D^{n+1}(t) + C E^{n+1}(t) D^n(t).
	\end{align*}
Denote $\mathcal{E}_n(t)= E^n(t) + \int_0^t D^n(s)\d s$ with $\mathcal{E}_n(0) = E^0(t) \equiv E(0) \le \tfrac{1}{2} M_0$. Then it follows that
	\begin{align}
	  & \tfrac{1}{2} E^{n+1}(s) \d s + \left( \tfrac{3}{5} - C M_0^{\frac{1}{2}} \right) \int_0^t D^{n+1}(s) \d s \\\no
	  \le\ & \tfrac{1}{2} E^{n+1}(0) + C t M_0 + C (M_0^{\frac{1}{2}} + 1) \int_0^t E^{n+1}(s) \d s
	  			 + C \int_0^t E^{n+1}(s) (D^n)^{\frac{1}{2}}(s) \d s \\\no
	  		 & + C \int_0^t (E^{n+1})^{\frac{1}{2}}(s) (D^n)^{\frac{1}{2}}(s) \d s
	  		   + C \int_0^t E^{n+1}(s) D^n(s) \d s.
	\end{align}
By the H\"older inequality again, we can infer that
	\begin{align*}
	  \int_0^t E^{n+1}(s) (D^n)^{\frac{1}{2}}(s) \d s
	  	\le\ & \int_0^t E^{n+1}(s) \d s + \int_0^t E^{n+1}(s) D^n(s) \d s \\
	  	\le\ & t \SE_{n+1}(s) + \SE_{n+1} \int_0^t D^n(s) \d s \\
	  	\le\ & t \SE_{n+1}(s) + \SE_n \SE_{n+1}(s),
	\end{align*}
and
	\begin{align*}
	  \int_0^t (E^{n+1})^{\frac{1}{2}}(s) (D^n)^{\frac{1}{2}}(s) \d s
	  	\le\ & \left\{ \int_0^t E^{n+1}(s) \d s \right\}^{\frac{1}{2}}
	  				 \left\{ \int_0^t D^n(s) \d s \right\}^{\frac{1}{2}}
	  	\le\ t^{\frac{1}{2}} \SE_{n+1}(s)^{\frac{1}{2}} \SE_n \\
	  	\le\ & \tfrac{1}{16 C} \SE_{n+1}(s) + C t \SE_n .
	\end{align*}
Similar process applied to the right-hand side terms in the above integral formulation will enables us to get
	\begin{align}
	  & \tfrac{1}{2} E^{n+1}(s) \d s + \left( \tfrac{3}{5} - C M_0^{\frac{1}{2}} \right) \int_0^t D^{n+1}(s) \d s \\\no
	  \le\ & \tfrac{1}{2} E^{n+1}(0) + C t M_0 + C t (M_0^{\frac{1}{2}} + 1) \SE_{n+1}(s) + C M_0 \SE_{n+1}(s) + \tfrac{1}{16} \SE_{n+1}(s) .
	\end{align}
Now we are ready to concludes that there exists some sufficiently small constant $M_0$ satisfying that $C M_0^{\frac{1}{2}} \le \tfrac{1}{10}$, and a time $T_*>0$ such that, for $ t \le T_*$,
	\begin{align}
	  ( \tfrac{7}{16} - C T_* (M_0^{\frac{1}{2}} + 1) - C M_0 ) \SE_{n+1}(s)
	  \le \tfrac{1}{4} M_0 + C T_* M_0,
	\end{align}
which leads us to the desired uniform-in-$n$ estimate
	\begin{align}
	  \sup_{t \in [0,T_*]} \E_{n+1} (t) \le M_0.
	\end{align}
Actually, $T_* \le M_0^{\frac{1}{2}} \le \tfrac{1}{32 C}$ suffices to match the requirement. This completes Lemma \ref{lemm:uniform esm}.

\subsection{Contraction in a low norm and local existence} % (fold)
\label{sub:contraction_in_low_norm}

% subsection contraction_in_low_norm (end)

Now we are left to prove the convergence. Set $\widetilde{\rho}^{n+1}= \rho^{n+1} -\rho^n, \widetilde{u}^{n+1}= u^{n+1} - u^n, \widetilde{f}^{n+1} = f^{n+1} -f^n $, then it follows from equation \eqref{eq:appro} that
	\begin{align}\label{eq:appro conv}
	  \begin{cases}
	  	\p_t \widetilde{\rho}^{n+1} + u^n \nabla_x \widetilde{\rho}^{n+1} + (1+\rho^n) \div \widetilde{u}^{n+1} = - \widetilde{u}^n \nabla_x \rho^n - \widetilde{\rho}^n \div u^n ,
		\\[5pt]
			\begin{aligned}
				\p_t \widetilde{u}^{n+1} & + u^n \cdot \nabla_x \widetilde{u}^{n+1} + \tfrac{P'(1+\rho^n)}{1+\rho^n} \nabla_x \widetilde{\rho}^{n+1}
				  + (\tfrac{P'(1+\rho^n)}{1+\rho^n} - \tfrac{P'(1+\rho^{n-1})}{1+\rho^{n-1}}) \nabla_x \rho^n
			\\[3pt]
				& = \tfrac{1}{1+\rho^n} \div \Sigma(\widetilde{u}^{n+1}) + \tfrac{1}{1+\rho^n} \div \int_{\mathbb{R}^3} \nabla_q U \otimes q \widetilde{g}^{n+1} M \d q \\
				& \quad + (\tfrac{1}{1+\rho^n} - \tfrac{1}{1+\rho^{n-1}}) \div \Sigma(u^n) + (\tfrac{1}{1+\rho^n} - \tfrac{1}{1+\rho^{n-1}}) \div \int_{\mathbb{R}^3} \nabla_q U \otimes q g^n M \d q ,
			\end{aligned}
	  \\[1em]
			\begin{aligned}
				\p_t \widetilde{g}^{n+1} & + u^n \cdot \nabla_x \widetilde{g}^{n+1} + \nabla_x u^n q \nabla_q \widetilde{g}^{n+1} + \mathcal{L} \widetilde{g}^{n+1}
				+ \widetilde{u}^n \cdot \nabla_x g^n + \nabla_x \widetilde{u}^n q \nabla_q g^n
			\\[3pt]
				& = - 2 \widetilde{g}^{n+1} \div u^n + \widetilde{g}^{n+1} \nabla_x u^n q \nabla_q U
				    - 2 (1+g^n) \div \widetilde{u}^n + (1+g^n) \nabla_x \widetilde{u}^n q \nabla_q U .
			\end{aligned}	
	  \end{cases}
	\end{align}	

Performing the $L^2$ energy estimate (with respect to $x$ or both $x$ and $q$) yields similar results as before, as follows:
	\begin{align}\label{esm:rho tilde}
	  \tfrac{1}{2} \tfrac{\d}{\d t} |\wr^{n+1}|^2_\lx
	    \ls\ & |\nabla_x u^n|_{L^\infty_x} |\wr^{n+1}|^2_\lx
	    		   + (1+ |\rho^n|_{L^\infty_x}) |\nabla_x \wu^{n+1}|_\lx |\wr^{n+1}|_\lx \\\no
	    		 & + |\nabla_x \rho^n|_{L^\infty_x} |\wu^n|_\lx |\wr^{n+1}|_\lx
	    		   + |\nabla_x u^n|_{L^\infty_x} |\wr^n|_\lx |\wr^{n+1}|_\lx \\\no
	    \ls\ & M_0^{\frac{1}{2}} |\wr^{n+1}|^2_\lx
	    			 + |\wr^{n+1}|_\lx |\nabla_x \wu^{n+1}|_\lx
	    			 + M_0^{\frac{1}{2}} |\wr^{n+1}|_\lx (|\wr^n|_\lx + |\wu^n|_\lx),
	\end{align}
and
	\begin{align}\label{esm:u tilde}
	  & \tfrac{1}{2} \tfrac{\d}{\d t} |\wu^{n+1}|_\lx^2
	  			+ (\mu |\nabla_x \wu^{n+1}|_\lx^2 + (\mu+\xi) |\div \wu^{n+1}|_\lx^2) \\\no
	  	\ls\ & |\nabla_x u^n|_{L^\infty_x} |\wu^{n+1}|_\lx^2
	  					+ |\wr^{n+1}|_\lx |\div \wu^{n+1}|_\lx
	  					+ |\nabla_x \rho^n|_{L^\infty_x} |\wu^{n+1}|_\lx (|\wr^{n+1}|_\lx + |\wr^n|_\lx) \\\no
	  			 & + |\div \Sigma(u^n)|_\lx |\wu^{n+1}|_\lx
	  			 	 + |\nabla_x \wu^{n+1}|_\lx \nm{\agl{q} \wg^{n+1}}_\lxq
	  			 	 + \nm{\agl{q} \nabla_x g^n}_\lxq |\wu^{n+1}|_\lx \\\no
	  	\ls\ & M_0^{\frac{1}{2}} (|\wu^{n+1}|_\lx^2 + |\wu^{n+1}|_\lx)
	  				 + |\wr^{n+1}|_\lx |\nabla_x \wu^{n+1}|_\lx
	  				 + \nm{\agl{q} \wg^{n+1}}_\lxq |\nabla_x \wu^{n+1}|_\lx \\\no
	  			 & + M_0^{\frac{1}{2}} |\wr^{n+1}|_\lx |\wu^{n+1}|_\lx,
	\end{align}
where we have used in the above two estimates the uniform bound $\sup_{t \in [0,T_*]} \E_{n+1} (t) \le M_0 <1$ and the Sobolev embedding inequality. At the same time, we can also deduce that
	\begin{align}\label{esm:g tilde}
	  & \tfrac{1}{2} \tfrac{\d}{\d t} \nm{\agl{q} \wg^{n+1}}_\lxq^2
	  			+ \nm{\agl{q} \nabla_q \wg^{n+1}}_\lxq^2 \\\no
	  	\ls\ & |\nabla_x u^n|_{L^\infty_x} \nm{\agl{q} \wg^{n+1}}_\lxq^2
	  	      + |\nabla_x u^n|_{L^\infty_x} \nm{\agl{q} \nabla_q \wg^{n+1}}_\lxq^2 \\\no
	  	    & + |\wu^n|_\lx \nm{\agl{q} \nabla_x g^n}_{L^\infty_x L^2_{Mq}} \nm{\agl{q} \wg^{n+1}}_\lxq
	  	     + |\nabla_x \wu^n|_\lx \nm{\agl{q}^2 \nabla_q g^n}_{L^\infty_x L^2_{Mq}} \nm{\agl{q} \wg^{n+1}}_\lxq \\\no
	  	    & + |\div u^n|_{L^\infty_x} \nm{\agl{q} \wg^{n+1}}_\lxq^2
	  	      + |\div \wu^n|_\lx \nm{\agl{q} \wg^{n+1}}_\lxq \\\no
	  	    & + |\div \wu^n|_\lx \nm{\agl{q} g^n}_{L^\infty_x L^2_{Mq}} \nm{\agl{q} \wg^{n+1}}_\lxq
	  	      + |\nabla_x \wu^n|_\lx \nm{\agl{q} \wg^{n+1}}_\lxq \\\no
	  	    & + |\nabla_x \wu^n|_\lx \nm{\agl{q} \nabla_q g^n}_{L^\infty_x L^2_{Mq}} \nm{\agl{q} \nabla_q \wg^{n+1}}_\lxq \\\no
	  	\ls\ & M_0^{\frac{1}{2}} (\nm{\agl{q} \wg^{n+1}}_\lxq^2 + \nm{\agl{q} \nabla_q \wg^{n+1}}_\lxq^2)
	  			   + M_0^{\frac{1}{2}} |\wu^n|_\lx \nm{\agl{q} \wg^{n+1}}_\lxq \\\no
	  			 & + \nm{\agl{q} \nabla_q g^n}_\hxqm{3} |\nabla_x \wu^n|_\lx \nm{\agl{q} \wg^{n+1}}_\lxq
	  			 	 + (1+M_0^{\frac{1}{2}}) |\nabla_x \wu^n|_\lx \nm{\agl{q} \wg^{n+1}}_\lxq.
	\end{align}

Let $\widetilde{E}^n(t)= E(\wr^n,\wu^n,\wg^n)(t)$, $\widetilde{D}^n(t)= D(\wr^n,\wu^n,\wg^n)(t)$, and $\widetilde{\mathcal{E}}_n(t)= \widetilde{E}^n(t) + \int_0^t \widetilde{D}^n(s)\d s$, then it follows from combining the above estimates \eqref{esm:rho tilde}-\eqref{esm:g tilde} that
	\begin{align}\label{esm:total tilde}
	  & \tfrac{1}{2} \tfrac{\d}{\d t} \WE^{n+1}(t) + \WD^{n+1}(t) \\\no
	  \le\ & C M_0^{\frac{1}{2}} [\WE^{n+1}(t) + \WD^{n+1}(t)]
	  			+ C M_0^{\frac{1}{2}} (\WE^{n+1})^{\frac{1}{2}} (t) [ (\WE^n)^{\frac{1}{2}} (t) + (\WD^n)^{\frac{1}{2}} (t) ] \\\no
	  		 & + C (\WE^{n+1})^{\frac{1}{2}} (t) (\WD^n)^{\frac{1}{2}} (t) [(D^n)^{\frac{1}{2}} (t) +1]
	  		   + C (\WE^{n+1})^{\frac{1}{2}} (t) (\WD^{n+1})^{\frac{1}{2}} (t) \\\no
	  		 & + C M_0^{\frac{1}{2}} (\WD^n)^{\frac{1}{2}} (t) (\WD^{n+1})^{\frac{1}{2}} (t).
	\end{align}
Notice that the last two terms can be controlled by the H\"older inequality, i.e.
	\begin{align*}
		(\WE^{n+1})^{\frac{1}{2}} (\WD^{n+1})^{\frac{1}{2}}
				\le\ & \tfrac{1}{4} \WD^{n+1} + C \WE^{n+1}, \\
	  M_0^{\frac{1}{2}} (\WD^n)^{\frac{1}{2}} (\WD^{n+1})^{\frac{1}{2}}
	  		\ls\ & M_0^{\frac{1}{2}} (\WD^n + \WD^{n+1}),
	\end{align*}
then we can get
	\begin{align}
    & \tfrac{1}{2} \tfrac{\d}{\d t} \WE^{n+1}(t)
    		+ \left( \tfrac{3}{4} - 2C M_0^{\frac{1}{2}} \right) \WD^{n+1}(t) \\\no
    	\le\ & C \WE^{n+1}(t) + C M_0^{\frac{1}{2}} \WD^n (t) + C M_0^{\frac{1}{2}} (\WE^{n+1})^{\frac{1}{2}} (t) (\WE^n)^{\frac{1}{2}} (t) \\\no
	  		 & + C (\WE^{n+1})^{\frac{1}{2}} (t) (\WD^n)^{\frac{1}{2}} (t) [(D^n)^{\frac{1}{2}} (t) +1].
	\end{align}
Then similar process as that of justifying the uniform boundedness of the approximate sequence will lead us to the integral formulation, by noticing the fact $\WE(0) \equiv 0$
	\begin{align}
	  \left( \tfrac{3}{4} - CT_* - C M_0^{\frac{1}{2}} \right)
	  	\sup_{t\in[0, T_*]} \WWE_{n+1}(t)
	  \le (C T_* + C M_0^{\frac{1}{2}} + C M_0^{\frac{1}{2}} T_*) \sup_{t\in[0, T_*]} \WWE_n(t).
	\end{align}
Consequently, the same choose of $M_0$ and $T_*$ as before yield the contraction property, i.e.
	\begin{align}
	  \sup_{t\in[0, T_*]} \WWE_{n+1}(t) \le c \sup_{t\in[0, T_*]} \WWE_n(t),
	\end{align}
with some constant $c<1$. Then combining with a standard compactness argument, we can get finally a unique local solution of the compressible polymeric system \eqref{eq:C-polymeric-pert}, which satisfies the bound $\sup_{t\in[0, T_*]} \E(t) \le M_0$ with initial datum $\E(0) \le \tfrac{1}{2} M_0$.

\subsection{Global existence by continuum argument} % (fold)
\label{sub:global_existence_by_continuum_argument}

% subsection global_existence_by_continuum_argument (end)

As the end step, aiming at proving the global-in-time solutions of the compressible polymeric system \eqref{eq:C-polymeric-pert}, we use one standard continuum argument. For that, we consider the initial datum $E(0) \le \tfrac{1}{2} M_0$ with $M_0 = \tfrac{1}{C_0} \min\{\ep,\, M_0\}$. Then we immediately get from Proposition \ref{prop:local exis} the unique local solution result on the time interval $t \in [0, T_*]$ with $T_*>0$, satisfying $\sup_{t\in[0, T_*]} \E(t) \le M_0$. For this local solution, the a priori estimate stated in Proposition \ref{prop:A-priori} yields that $E(T_*) \le C_0 E(0) \le \tfrac{1}{2} C_0 M_0$, using Proposition \ref{prop:local exis} again, it follows that the time interval can be extended by another $T_*$, and we have $\sup_{t\in[T_*, 2T_*]} \E(t) \le C_0 M_0$. So it holds that
	\begin{align*}
	  \sup_{t\in[0, 2T_*]} \E(t) \le C_0 M_0.
	\end{align*}
By using the a priori estimate in Proposition \ref{prop:A-priori} again, we get $E(2T_*) \le C_0 E(0)$. This will lead to the existence result on time interval $t \in [0, 2 T_*]$. Repeating this bootstrap argument, we are led to the global existence of classic solution of the compressible polymeric system \eqref{eq:C-polymeric-pert} near equilibrium, moreover, it satisfies the energy bound $\sup_{t\in[0, +\infty)} \E(t) \le C_0 M_0 < \ep$.

%% ------------------------------------------------------------------------ %%

%% ------------------------------------------------------------------------ %%

% \bigskip
% \phantomsection
% \addcontentsline{toc}{section}{\refname}

%\bibliographystyle{unsrtnat}
%\nocite{*}
%\bibliography{reference}

\begin{thebibliography}{99}

\bibitem{BAH-1987-BOOK} R.~B. Bird, R.~Amstrong and O.~Hassager.
\newblock {\em Dynamics of polymeric liquids Vol. 1,}.
\newblock Wiley, New York, 1977.

\bibitem{BCAH-1987-BOOK} R.~B. Bird, C.~Curtiss, R.~Amstrong and O.~Hassager.
\newblock {\em Dynamics of polymeric liquids, Kinetic Theory Vol. 2,}.
\newblock Wiley, New York, 1987.

\bibitem{CG-2006-CPDE} J. A. Carrillo and T. Goudon.
\textit{Stability and asymptotic analysis of a fluid-particle interaction model}. {Comm. Partial Differential Equations} {\bf 31} (2006), no. 9, 1349-1379.

\bibitem{CM-2001-SIMA} J.-Y. Chemin and N.~Masmoudi.
\newblock \textit{About lifespan of regular solutions of equations related to viscoelastic fluids}.
\newblock {SIAM J. Math. Anal.} {\bf 33} (2001), no. 1, 84-112.

\bibitem{CFTZ-2007-CMP} P. Constantin, C. Fefferman, E. S. Titi and A. Zarnescu.
\textit{Regularity of coupled two-dimensional nonlinear Fokker-Planck and Navier-Stokes systems}. {Comm. Math. Phys.} {\bf 270} (2007), no. 3, 789-811.

\bibitem{CM-2008-CMP} P. Constantin and N. Masmoudi.
\textit{Global well-posedness for a Smoluchowski equation coupled with Navier-Stokes equations in 2D}. {Comm. Math. Phys.} {\bf 278} (2008), no. 1, 179-191.

\bibitem{ELZ-2004-CMP} W. E, T. Li and P. Zhang. \textit{Well-posedness for the dumbbell model of polymeric fluids}. {Commun. Math. Phys.} {\bf 248} (2004), no. 2, 409-427.

\bibitem{HKL-2010-DCDS} Y. Hyon, D. Y. Kwak and C. Liu. \textit{Energetic variational approach in complex fluids: maximum dissipation principle}. {Discrete Contin. Dyn. Syst} {\bf 26} (2010), no. 4, 1291-1304.

\bibitem{JLL-2004-JFA} B.~Jourdain, T.~Leli{\`e}vre and C.~Le~Bris.
\newblock \textit{Existence of solution for a micro-macro model of polymeric fluid: the
  {FENE} model}.
\newblock {J. Funct. Anal.} {\bf 209} (2004), no. 1, 162-193.

\bibitem{JLL-2016-Deborah} N. Jiang, F-H, Lin and Y.L. Luo. \textit{Zero Deborah number limit of the micro-macro polymeric flows}. Preprint(2016).

\bibitem{LLZ-2008-ARMA} Z.~Lei, C.~Liu and Y.~Zhou.
\newblock \textit{Global solutions for incompressible viscoelastic fluids}.
\newblock {Arch. Ration. Mech. Anal.} {\bf 188} (2008), no. 3, 371-398.


\bibitem{LLZ-2007-CPAM} F.H. Lin, C. Liu and P. Zhang. \textit{On the Micro-Macro Model for Polymeric Fluids near Equilibrium}. {Comm. Pure Appl. Math.} {\bf 60} (2007), no. 6, 838-866.

\bibitem{LZZ-2008-CMP} F.-H. Lin, P.~Zhang and Z.~Zhang.
\newblock \textit{On the global existence of smooth solution to the 2-{D} {FENE}
  dumbbell model}.
\newblock {Comm. Math. Phys.} {\bf 277} (2008), no. 2, 531-553.

\bibitem{LM-2000-CAM} P.-L. Lions and N.~Masmoudi.
\newblock \textit{Global solutions for some {O}ldroyd models of non-{N}ewtonian flows}.
\newblock {Chinese Ann. Math. Ser. B} {\bf 21} (2000), no. 2, 131-146.

\bibitem{Liu-2009-NOTES} C. Liu.
\textit{An Introduction of Elastic Complex Fluids: An Energetic Variational Approach}. Multi-Scale Phenomena In Complex Fluids: Modeling, Analysis and Numerical Simulation. 2009, 286-337.

\bibitem{Majda-1984-BOOK} A.~J. Majda. {Compressible Fluid Flow and Systems of Conservation Laws in Several Space Variables}. Vol. 53. Springer Science \& Business Media, 1984.

\bibitem{MB-2002-BOOK} A.~J. Majda and A.~L. Bertozzi. {Vorticity and incompressible flow}. Vol.~27 of Cambridge Texts in Applied Mathematics, Cambridge University Press, Cambridge, 2002.

\bibitem{Masmdi-2008-CPAM} N.~Masmoudi.
\newblock \textit{Well-posedness for the {FENE} dumbbell model of polymeric flows}.
\newblock {Comm. Pure Appl. Math.} {\bf 61} (2008), no. 12, 1685-1714.

\bibitem{Masmdi-2013-Invent} N.~Masmoudi. \textit{Global existence of weak solutions to the FENE dumbbell model of polymeric flows}. {Invent. Math.} {\bf 191} (2013), no. 2, 427-500.

\bibitem{MN-1983-CMP} A. Matsumura and T. Nishida. \textit{Initial-boundary value problems for the equations of motion of compressible viscous and heat-conductive fluids}. {Comm. Math. Phys.} {\bf 89} (1983), 445-464.

\bibitem{MV-2006-CPDE} A. Mellet and A. Vasseur.
\textit{Global weak solutions for a Vlasov–Fokker–Planck/compressible Navier-Stokes system of equations}.{Math. Models Methods Appl. Sci.} {\bf 17} (2007), no. 7, 1039-1063.

\bibitem{Moser-1966-1} J. Moser. \textit{A rapidly convergent iteration method and non-linear partial differential equations-I}. {Annali della Scuola Normale Superiore di Pisa-Classe di Scienze}, {\bf 20} 1966, no. 2, 265-315.

\bibitem{OT-2008-CMP} F. Otto and A. E. Tzavaras.
\textit{Continuity of Velocity Gradients in Suspensions of Rod–like Molecules}. {Comm. Math. Phys.} {\bf 277} (2008), no. 3, 179-191.

\bibitem{SL-2009-DCDS} H. Sun and C. Liu. \textit{On energetic variational approaches in modeling the nematic liquid crystal flows}. {Discrete Contin. Dyn. Syst} {\bf 23} (2009), no. 1-2, 455-475.

\bibitem{TM-2011-BOOK} E. B. Tadmor and R. E. Miller. Modeling materials. Continuum, atomistic and multiscale techniques, Cambridge: Cambridge University Press, 2011.

\bibitem{Taylor-2011-BOOK} M.~E. Taylor. Partial differential equations {III}. {N}onlinear equations, vol. 117 of Applied Mathematical Sciences, 2nd ed., Springer, New York, 2011.

\bibitem{ZZ-2006-ARMA} H.~Zhang and P.~Zhang.
\newblock \textit{Local existence for the {FENE}-dumbbell model of polymeric fluids}.
\newblock {Arch. Ration. Mech. Anal.} {\bf 181} (2006), no. 2, 373-400.

\end{thebibliography}

\end{document}